\documentclass[11pt]{article}
\usepackage{amsmath}
\usepackage{amssymb}
\usepackage{amsthm}
\usepackage{anysize}

\usepackage{mathrsfs}
\usepackage{natbib}
\usepackage{xspace}

\usepackage[createShortEnv]{proof-at-the-end}

\usepackage{dsfont}
\usepackage{enumerate}
\usepackage{graphicx}
\usepackage{floatrow}
\usepackage{bm}
\usepackage{multicol}
\usepackage{multirow}
\usepackage[dvips]{epsfig}

\usepackage{float}
\usepackage{color}
\usepackage[latin1]{inputenc}
\usepackage[english]{babel}
\usepackage{xcolor}

\usepackage{bbm}
\usepackage{mathtools}

\numberwithin{equation}{section}

\theoremstyle{plain}

\newtheorem{cor}{Corollary}

\theoremstyle{definition}
\newtheorem{defn}{Definition}
\theoremstyle{definition}

\newtheorem{esmp}{Example}

\theoremstyle{remark}
\newtheorem{oss}{Remark}

\marginsize{25mm}{25mm}{25mm}{25mm}

\newcommand{\I}{\mathds{1}}

\newcommand{\R}{\mathbb R}

\newcommand{\Cov}{\mathrm{Cov}}

\renewcommand \|{|}

\begin{document}

\title{Thorin processes and their subordination\footnote{Part of this paper has been presented at the 2024 Meeting in Probability and Mathematical Statistics in Rome. The author is grateful to Nikolay Leonenko and other participants for the useful comments. Any error is this author's own.} \footnote{This document has been redact without any kind of generative artificial intelligence or large language models assistance, in either content or form.}
}
\author{\\ Lorenzo Torricelli\footnote{University of Bologna, Department of Statistical Sciences ``P. Fortunati''. Email: lorenzo.torricelli2@unibo.it}   }
\date{\today}

\maketitle

\begin{abstract} 

A Thorin process is a stochastic process with independent and stationary increments whose  laws are weak limits of finite convolutions of gamma distributions. Many popular L\'evy processes fall under this class. The Thorin class can be characterized by a representing triplet that conveys more information on the process compared to the L\'evy triplet. In this paper we investigate some relationships between the Thorin structure and the process  properties, and find that the support of the Thorin measure characterizes the existence of the critical exponential moment,  as well as the asymptotic equivalence between the L\'evy tail function and the complementary distribution function. Furthermore, it is illustrated how univariate Brownian subordination  with respect to Thorin subordinators produces Thorin processes whose representing  measure is given by a pushforward with respect to a hyperbolic function, leading to arguably easier formulae compared to the Bochner  integral determining the L\'evy measure. We provide a full account of the theory of multivariate Thorin processes, starting from the Thorin--Bondesson representation for the characteristic exponent, and highlight the roles of the Thorin measure in the analysis of density functions, moments, path variation and  subordination. Various old and new  examples are discussed. We finally detail a treatment of subordination of gamma processes with respect to negative binomial subordinators which is made possible by the Thorin--Bondesson representation.
\end{abstract}


\noindent {\bf{AMS}}: 60G51, 60E07


\section{Introduction} 

 An infinitely divisible (i.d.) distribution  is one that writes as an $n$-fold convolution of some other probability distribution for any arbitrary $n \in \mathbb N$. The concept originates from the work of de Finetti in the late 1920s, and was first rigorously investigated in \citet{lev:34}. It is hard to overstate the importance of i.d. distributions in applications to physical and social sciences. Any stochastic process with stationary and independent increments, all-important properties for statistical inference, must possess this property. Conversely, a distribution of this type naturally gives raise to one such process, whose class is named after P. L\'evy. The popularity of L\'evy processes for mathematical modeling is firmly established. In actuarial mathematics for example, the most basic model for claim arrivals and assessing ruin probabilities is the compound Poisson process, that, together with the Brownian motion, represents the most basic example of L\'evy process. In financial economics L\'evy processes were popularized starting in the 90s (\cite{mad+sen:90}, \cite{mad+al:98}) as  models for asset returns, capturing non-normality and the discontinuous trajectories those might exhibit. Small time central limit theorems classically find application within this class, with Gaussian or stable limiting regimes attained depending upon the existence of variance of the sample distribution. From an analytical standpoint, the characteristic function (c.f.) of i.d. probability laws is canonically obtained from a triplet consisting of a vector, a symmetric nonnegative definite matrix and a real measure (the ``L\'evy measure''), through the so-called L\'evy-Khincthine representation. In particular the L\'evy measure is an analytical determinant of many distributional and path properties of the process, alleviating the typical problem of the absence of closed formulae for the distribution functions themselves. Systematic investigations in the theory of i.d. probability distributions has been an active area of research throughout all the second half of the 20th century, culminating around the turn of the century in classic monographs such as 	
 \citet{ber:96}, \citet{sat:99}, \citet{ste+vh:03} and \citet{app:09}. 

Between the 60s and 70s a specific i.d. subclass came to prominence,  the $L$-class, or that of the self-decomposable  distributions. 
The L\'evy measure of self-decomposable distributions comes in a specific absolutely continuous form, incorporating a radially decreasing function conveying more information than the L\'evy measure itself. The implications that  self-decomposability bears on the underlying probability law were studied among the others by \citet{zol:63}, \citet{wol:71}, \citet{yam:78}, \citet{sat+yam:78}. One particularly important aspect is absolute continuity, addressed in \citet{sat:80}, \citet{sat:82}, \citet{yam:83}. Also, such class appears naturally in the study of Ornstein-Uhlenbeck processes and their stationary log-transform, as in \citet{wol:82}, \citet{jur+ver:83}, \citet{bn+she:02}.

A strict subclass of the univariate self-decomposable class can be obtained by looking at the smallest set containing all of the positive and real supported gamma  distributions, their finite convolutions and weak limits. Such class was introduced by Olof Thorin  in  1977-1978 en route to his proof of the infinite-divisibility of the lognormal distribution (\citet{tho:77}), under the name of (extended) generalized gamma  convolutions ((E)GGC)
Many well-known  distributions belong to this class, including finite gamma convolutions, stable and tempered stable laws (\citet{kop:95}), logistic distributions, Gaussian, (generalized) inverse Gaussian distributions, and many others.   For a thorough introduction to univariate Thorin laws we refer to the definitive monograph \citet{bon:12}; a chapter on GGC laws is also contained in \citet{ste+vh:03}. As far as this author is aware, multivariate EGGC distributions were first considered in \citet{bns+al:06}, in relation with their representation as laws of certain classes of stochastic integrals. 
In \citet{gri:08} and \citet{gri:11} the author studies possible extensions of the positive Thorin laws by letting the weak limits being taken in the Tweedie/positive tempered stable class, instead of among gamma distributions. For applications,  \citet{buc+al:17} suggest modeling financial securities by means of  subordination of multivariate Brownian motions with respect to multivariate GGCs (i.e. Thorin processes supported on the positive orthant).  An interesting survey of some distributional properties enjoyed by GGCs is provided in 
\citet{jam+al:08}.   \citet{bur:14}  prove an invariance property of the GGC class with respect to a transformation of processes based on Kendall's identity.
A  recent discussion on the statistical estimation of multivariate generalized gamma convolutions using Laguerre basis expansion is \citet{lev+al:21}.

The key feature of the Thorin laws is that when written in polar coordinates their  canonical function is completely monotone along each direction, that is, by the famous theorem of N. Bernstein, they are radially the Laplace transform of a positively supported  -- not necessarily finite -- measure. Such measure, named after Thorin himself, is a representing measure for the c.f. of the law, and a Thorin characteristic triplet for the c.f. paralleling that of L\'evy can be identified.  In  applied sciences, when modeling phenomena with stationary and independent innovations,  Thorin processes appear to be more of a rule than an exception.   The view we take throughout this paper, is that even if the L\'evy measure can be worked out by Laplace-transforming the Thorin measure, it may be at times convenient to study the process from a pure Thorin point of view. A circumstance in which this happens is  when subordinating Brownian motions to an independent tempered geometric stable subordinator (\citet{tor+al:21}), whose Thorin measure is readily found, while the corresponding L\'evy measure is not known. Compare with  Example \ref{es:NL} below.

A representation Theorem in terms of the Thorin triplet in the univariate case has been given as early as \citet{tho:78} and \citet{bon:12}. \citet{gri:08} also describes the c.f. for his generalized multivariate Thorin laws, and in principle a full multidimensional representation for the c.f.s of Thorin distributions should be obtained as a specialization of those. It is however hard (at least for this author) to reconcile  the one-dimensional case of such general expression with what was originally established by \citet{tho:78} and \citet{bon:12}.  We therefore provide in Section 2 a  full representation Theorem for multivariate Thorin c.f.s in terms of the Thorin triplet,  christened here ``Thorin--Bondesson representation'', that transparently reconciles with the univariate case.  We subsequently introduce the concept of  L\'evy-Thorin processes, as those processes with stationary and independent increments having unit time marginals in the Thorin class, which constitutes the main focus of this paper.
 
In Section 3 we reformulate some distributional and path variation properties of Thorin processes in terms of their Thorin, as opposed to L\'evy, triplet. We then exhibit  regularity results on the density functions that can be attained by mirroring those of \citet{sat+yam:78} for self-decomposable distributions. The emphasis is here on the fact that a zero Blumenthal-Getoor (BG) index affects both the order of differentiability and the boundedness around zero of the density. An important finding is that the support of the Thorin measure relates to the existence of the critical exponential moment for the distribution. In turn, this plays a primary role in the large-value asymptotics of the distribution function, since  the existence of the moment generating function at the critical exponents makes it possible to apply an important Theorem (\citet{wat:08}, \citet{sgi:90}, \citet{pak:04}) guaranteeing that the tails of the L\'evy measure and those of the complementary distribution function are equivalent.

In Section 4 we revisit some important instances of Thorin distributions, determine their Thorin measure, and see how they fit in in the theory exposed.  We focus on \emph{(a)} finite Gamma convolutions, \emph{(b)} stable and \emph{(c)} tempered stable (including geometric stable)  laws,  \emph{(d)} generalized-logistic type laws.

Section 5 is the core of the paper and deals with Thorin subordination. In a one-dimensional set up we prove a series of results characterizing Brownian subordination with respect to GGC processes. The most important insights are that on the one hand, in distribution, a drifted Brownian motion subordinated to a GGC is always a Thorin process. On the other, as a converse, a Thorin process with shift-symmetric Thorin measure is distributionally equal to a GGC-subordinated Brownian motion with drift. Another advantage of using the Thorin representation of c.f.s is that it paves the way for a natural subordination theory of gamma and bilateral gamma processes,  producing  Thorin processes with zero BG index, provided that the subordinator is chosen as a finite activity L\'evy process with lattice-valued negative-binomial distribution (and weak limits of convolutions thereof). The proofs are found in the last section. 

\bigskip

\section{Thorin distributions and processes}\label{sec:Thordists}

As customary $\R^d$ is the $d$-dimensional Euclidean space, and with $\R_+^d$ we indicate  the closed positive orthant. A Borel measure (or simply measure) $\eta$ on $D \subseteq \R^d$ is a measure on the measurable space $(D, \mathcal B(D))$ where $\mathcal B(D)$ is the Borel sigma-algebra on $D$. Integration of an integrable function $f$ with respect to one such measure is understood in the Lebesgue sense, unless otherwise stated, and for $a,b \in D \subseteq \R$, $f: D \subseteq \mathbb R\rightarrow \R$, $\int_a^b f(x) \eta(dx) $ is a shorthand for $\int_{[a,b]} f(x) \eta(dx)$.  If $d=1$ we let $F_\eta(x)=(-\infty, x]$, $\overline \eta(x)=(x, \infty)$ respectively the distribution and tail functions of $\eta$. If $\eta$ is absolutely continuous (with respect to the Lebesgue measure) with abuse of notation we denote its density by $\eta(x)$, $x \in D \subseteq \R^d$
For a Borel measurable (or simply measurable) function $g:D \subseteq \R^d \rightarrow \R^d$ and a measure $\eta$ on $\R^d$, $g_\sharp \eta$ indicates the pushforward measure given by $g_\sharp \eta(B)=\eta(g^{-1}(B))$, $B \in  \mathcal B(\R^d)$. If $\theta \in \R^d$, then $s^\theta:\R^d \rightarrow \R^d$ is by definition the function $s^\theta(x)=x-\theta$. Hence $s^\theta_\sharp \eta$ i.e. $s^\theta_\sharp \eta(B)=\eta(B+\theta)$ (Minkowksi sum), is the left translation of $\eta$ by $\theta$ and is indicated by $\eta^\theta$. For a measurable function $f:D \subseteq \R^d \rightarrow \R^d$ we let $\check f$  be $\check f(x)=f(-x)$. Similarly, for a measure $\eta$ on $\R^d$ its reflection is the measure $\check \eta(B)=\eta(-B)=i_\sharp \eta(B)$, $i(x)=-x$. Given two measures $\eta_1$ and $\eta_2$ their convolution is denoted $\eta_1 * \eta_2$. 

We use $\langle \cdot, \cdot \rangle$ for the standard scalar product on $\R^d$ and $|\cdot |$ for the Euclidean $d$-norm.
With distribution or law we mean throughout a probability distribution and probability law. For a law $\mu$ on $\mathbb R^d$ we let $\hat \mu(z)=\int_{-\infty}^\infty e^{i\langle z, x\rangle}\mu(dx)$, $z \in \mathbb R^d$, be its c.f.. The Laplace transform, wherever defined in $D \subseteq \mathbb C^d$  is  indicated $\mathscr L(\mu, s)$, $s \in D \subseteq \mathbb C^d$. The law of a random variable (r.v.) $X$ is denoted by $\mathcal L(X)$. 
  With $X \sim \mu$ and $X$ a r.v. and $\mu$ a probability law, we mean ``distributed as''. For two r.v.s $X$, $Z$, $X \sim Z$ means instead identity in distribution. For a measure $\eta$ its support is written supp$\, \eta$.  The symbol $\delta_a$ with $a \in \R^d$ stands for the Dirac delta distribution on $\R^d$ concentrated in $a$.  Notice that $\eta^\theta=\eta *\delta_{-\theta}$, because convolution of measures is the pushforward of the product measure with respect to the addition map $\pmb +: \R^d \times \R^d \rightarrow \R^d$, $\pmb{+}(x,y)=x+y$, i.e.  $\eta *\delta_{-\theta}= \pmb +_\sharp(\eta \times \delta_{-\theta})(B)=\eta(B+\theta)$, for all $B \in \mathcal B(\R^d)$. 
 Classes of distributions are denoted with $C(S)$, $C$ the class denomination, and $S$ some subset of the $d$-Euclidean space, typically $S=\R^d, \R^d_+, \R, \R_+$. 

A distribution $\mu$ on $\R^d$ is said to be infinitely divisible (i.d.)  if for all $z \in \mathbb R^d$
\begin{equation}\label{eq:LK}
\psi(z):=\log \hat \mu(z)=  \exp\left(i  \langle  z ,a \rangle   -  \frac{1}{2} \langle \Sigma z, z\rangle +\int_{\mathbb R^d} \left( e^{i \langle z,  x\rangle} -1 -   i  \langle  z,  x\rangle \I_{\{|  x | \leq 1\}}\right)\nu( d x)\right), 
\end{equation}
 with $\log$ the principal branch of the complex logarithm, $a \in \mathbb R^d$, $\Sigma$ a symmetric nonnegative definite square-matrix of order $d$, $\nu$  a measure on $\mathbb R^d$ such that $\nu(\{0\})=0$, $\int_{\mathbb R^d} (| x |^2 \wedge 1)\nu( dx) <\infty$. The function $\psi:\R^d \rightarrow \mathbb C$ is called the characteristic exponent. The class ID is that of the i.d. distributions.   For $\mu \in$ ID$(\R_+)$ it is convenient to use the Laplace exponent $\phi(s)=\log \mathscr L( \mu,s)=\psi(i s)$, $s  \geq 0$. To indicate the exponents $\psi$, $\phi$ as referred to a specific  r.v. $X$ or distribution $\mu$ we use subscripts i.e. $\psi_X$, $\phi_X$, $\psi_\mu$, $\phi_\mu$ .  
  The measure $\nu$ is called the \emph{L\'evy measure} and the triplet $(a, \Sigma, \nu)$ the \emph{L\'evy triplet}. 
 
One of the most important ID$(\R_+)$, respectively ID$(\R)$, distributions are the  gamma distributions, respectively, bilateral gamma distributions,  which are i.d. with triplet $(a,0, b e^{-c x}/x dx)$, $a \geq 0, b,c >0$, and  $\Big(a,0, (b_+ x^{-1}e^{-c_+ x}\I_{\{x  >0\}} + b_-|x|^{-1} e^{c_- |x|}\I_{\{x <0\}}) dx \Big),  a \in \R, b_\pm, c_\pm >0$ , respectively. If we need to emphasize the parameters we write $\Gamma(b,c)$ and Bil$\Gamma(b_+, c_+, b_-, c_-)$, ignoring the location value $a$.

A stochastic process $X=(X_t)_{t \geq 0}$ taking values in $\R^d$ is a \emph{ L\'evy process} on $\R^d$ if $X$ is stochastically continuous with stationary and independent increments.  Equivalently, $X=(X_t)_{t \geq 0}$ is a L\'evy process on $\R^d$ if and only if $\mathcal L (X_1) \in $ ID$(\R^d)$, and if this happens, $\mathcal L (X_t) \in$ ID$(\R^d)$ for all $t > 0$.  See \citet{sat:99} for details. We will omit mention of the range.

We introduce the multivariate Thorin class following \citet{bns+al:06}. However, in order to remain in line with the original spirit of Thorin and Bondesson, 
we interchange the definition and the characterization found therein (Definition 2.1 and Theorem F).

Denote by $\Gamma_e(\R^d)$ the class of \emph{elementary multivariate gamma distributions}, i.e. $\mu \in \Gamma_e(\R^d)$ if and only if $\mu  \sim V x$ where $V$ is a $\Gamma(\R_+)$ r.v. and $x \neq 0 \in \R^d$. 

\begin{defn}\label{eq:EGGCdef} We say that $\mu$ is a  \emph{$d$-dimensional extended generalized gamma convolution} (EGGC) or that $\mu$ belongs to the
 \emph{Thorin class}, and write  $\mu \in $ $T (\R^d)$, if there exists a sequence $(\mu_{n})_{n \geq 1}$, with $\mu_n=\mu_{1,n} * \ldots *\mu_{k_n,n}$, $ \mu_{i,n} \in \Gamma_e(\R^d)$, $k_n \in \mathbb N$, $i=1,\ldots, k_n$,  such that $\mu_n \rightarrow \mu$ weakly. 
\end{defn}
 
 In other words, $T(\R^d)$ is the smallest subclass in ID$(\R^d)$ containing all the elementary gamma distributions, their finite convolutions and weak limits.
 
 Let us focus on the $d=1$ case, and summarize the  known representation results. 
Following \citet{tho:78} and \citet{bon:12}, Section 7, an EGGC distribution $\mu$ on the real line, i.e. $\mu \in T(\R)$, is fully described by a triplet $(b, \theta, \tau)$ with $c \in \mathbb R$, a quadratic characteristic $\theta \geq 0$ and a measure $\tau$, such that  
\begin{equation}\label{eq:GGCLK}
\psi_\mu(z)=  i  z b + \frac{z^2}{2} \theta +\int_{\mathbb R}\left(\log\left(\frac{x}{x-iz}\right) - \frac{i z}{ x} \I_{\{| x | \geq 1\}}\right) \tau(d  x), \quad z \in \R
\end{equation}
provided that
\begin{equation}\label{eq:GGCtm}\begin{array}{lr}
 \displaystyle{\int_{\{|x| \leq 1\}} |\log|x||\tau(dx)<\infty}, & \displaystyle{\int_{\{|x| \geq 1\}} \frac{\tau(dx)}{x^2}} < \infty.
\end{array}
\end{equation}
Notice that compared to the mentioned references we are changing the truncation function in order to align it to the  modern choice of $x\I_{\{|x| \leq 1 \}}$ for the L\'evy measure truncation (compare with \eqref{eq:LK}).

The triplet $(b, \theta, \tau)$ is an example of Thorin triplet. The GGC subclass $T(\R_+) \subset T(\R)$  can instead be represented by the pair $(b,\tau)$, $b \geq 0$, where now the following Laplace representation holds
\begin{equation}\label{eq:thorinLK}
\phi_\mu(s)=   s b  +\int_0^\infty \log\left(\frac{x}{x + s}\right) \tau(d  x), \quad \mbox{Re}(s)>0
\end{equation}
provided that
\begin{equation}\begin{array}{lr}
 \displaystyle{\int_0^1 |\log x|\tau(dx)<\infty}, & \displaystyle{\int_1 ^\infty \frac{\tau(dx)}{x}} < \infty.
\end{array}
\end{equation}
In the following we extend these representations to the multivariate case.

 We need first to recall some results and definitions. To begin with, multivariate i.d. distributions are typically easier to study in polar coordinates. We have the following important proposition.
\begin{theoremEnd}[proof at the end,
	no link to proof]{prop}\label{prop:spherrep} Let $\mu \in $ {\upshape ID}$(\R^d)$, 
and denote $S^{d-1}=\{u \in \R^d  \, | \, \|u \|=1\}$.  There exists a finite measure $\lambda$ on $S^{d-1}$  and a family $\{m_u, \, u \in S^{d-1} \}$ of $\sigma$-finite measures on $\mathbb R_+$, respectively the ``spherical'' and ``radial'' components,  such that for all $B \in \mathcal B(\mathbb R_+\setminus \{ 0\})$, $m_u(B)$ is measurable as a function of $u$, and
\begin{equation}\label{eq:mupolar}
\nu(B)=\int_{S^{d-1}} \lambda(du) \int_{0}^\infty \I_{B}(r u )m_u(dr), \qquad \forall B \in \mathcal B(\mathbb R^d \setminus \{ 0\} ).
\end{equation}

Furthermore $\lambda$ and $m_u$ are uniquely determined up to multiplication by a positive  measurable function, i.e. for any $\lambda'$, $m_u'$ satisfying \eqref{eq:mupolar}, there exists a Borel-measurable function $c:S^{d-1} \rightarrow \R_+$ such that
\begin{align}
\lambda'(du)&=\lambda(du) c(u) \\
c(u) m_u'(dr)&=m_u(dr).
\end{align}
In particular, $\lambda$ can always be chosen as a probability measure.
\end{theoremEnd}

\begin{proofE}
 \citet{bns+al:06}, Lemma 2.1,   or \citet{gri:08}, Proposition 2.
\end{proofE}
Then if  $\mu \in$ ID$(\R^d)$ with L\'evy measure $\nu$, the pair $(\lambda, \{m_u , u \in S^{d-1} \})$ is a polar coordinate representation of $\nu$, uniquely determined up to a measurable function on $S^{d-1}$.  

 The  class SD of self-decomposable distributions can be conveniently introduced using polar coordinates.
\begin{defn}
We say that $\mu \in$ ID$(\R^d)$ is \emph{self-decomposable}, and write $\mu \in$ SD$(\R^d)$ if its radial component $\{ m_u, \, u \in S^{d-1} \}$ satisfies
\begin{equation}\label{eq:polarmease}
m_u(dr)=\frac{k(u,r)}{r} dr 
\end{equation}
 where $k: S^{d} \times \mathbb R_+ \rightarrow \R_+ $ is decreasing  and right-continuous in $r$ and measurable in $u$. The function $k$  is called the \emph{canonical function}.
\end{defn}
A L\'evy process $X=(X_t)_{t \geq 0}$ is then said to be self-decomposable if $\mathcal L(X_1) \in$ SD$(\R^d)$, or equivalently $\mathcal L(X_t) \in$ SD$(\R^d )$ for all $t >0$.

Recall that a \emph{completely monotone} (c.m.) function $f:\mathbb (0, \infty) \rightarrow \mathbb R$ is an infinitely-differentiable function such that
\begin{equation}
(-1)^{n}f^{(n)}(x) \geq 0,
\end{equation}
 with the parenthetical exponent indicating the $n$-th derivative.  A famous theorem by S.N. Bernstein states that a function is c.m. if and only if there exists a positive Borel measure $\eta$ such that 
\begin{equation}\label{eq:Bernsthm}
f(x)=\int_0^\infty e^{- s x }\eta(ds), \qquad x>0
\end{equation}
i.e.  $f$ is the Laplace transform of some such measure. The measure $\eta$ needs not to be finite, and consequently $f$ does not necessarily complete to a continuous function on $[0, \infty)$.

\begin{theoremEnd}[proof at the end,
	no link to proof]{thm}\label{thm:ThorCM}
A law $\mu$ is such that $\mu \in T(\R^d)$ if and only $\mu \in$ {\upshape SD}$(\R^d)$ and  its canonical function $k$ is c.m. in the variable $r$ for all $u \in S^{d-1}$, i.e.
\begin{equation}\label{eq:thorinm}
k(u,r)=\int_0^{\infty} e^{-s r} \tau_u(ds),  \qquad r>0
\end{equation}
for a family of positive measures $\{\tau_u, \, u \in S^{d-1}\}$ on $[0, \infty)$, such that $\tau_u(B)$ is measurable as a function of $u$ for all $B \in \mathcal B(\R_+)$. 
\end{theoremEnd}

\begin{proofE}
The claim follows from Definition \ref{eq:EGGCdef}  and   Definition 2.2 and Theorem F of \citet{bns+al:06}. 
\end{proofE}

Observe that as in \eqref{eq:mupolar} we can associate  to $\mu \in T(\R^d)$  the following Borel measure on $\R^d$
\begin{equation}\label{eq:ThorinM}
\tau(B)=\int_{S^{d-1}} \lambda(du) \int_{0}^\infty \I_{B}(s u )\tau_u(ds), \qquad \forall B \in \mathcal B(\mathbb R^d \setminus \{ 0 \}  ).
\end{equation}
We do not exclude $\tau=0$, which is the case if and only if  $k(u,r)=0$ i.e.  $\nu=0$ and $\mu$ is a multivariate Gaussian law.

In the case $d=1$ then $S^0=\{ 1,-1\}$ and $\lambda(du)=\lambda_+ \delta_1 + \lambda_+ \delta_{-1}$ with $\lambda_+, \lambda_- \geq 0$. In such case we rename $\tau_1=\tau_+$, $\tau_{-1}=\tau_-$ and absorb the costants $\lambda_\pm$ in $\tau_\pm$  obtaining  \begin{equation}\label{eq:tau1d}
\tau(B)=\tau_+(B \cap \R_+) +\check \tau_-(B \cap \R_-), \qquad B 
\in \mathcal B(\R \setminus \{ 0 \} )
\end{equation} which recovers the measure $\tau$ in  \eqref{eq:GGCLK} originally introduced by Thorin.

  We can now proceed with the main representation theorem.
\begin{theoremEnd}[proof at the end,
	no link to proof]{thm}[\textbf{Thorin--Bondesson representation}]\label{thm:olddef}
Let $\mu \in  T(\R^d)$ with L\'evy triplet $(a, \Sigma, \nu)$, where $\nu$ has polar coordinate rperesentation $(\lambda, \{m_u , u \in S^{d-1} \})$. Then the measure $\tau$ in \eqref{eq:ThorinM} satisfies \begin{equation}\label{eq:thorconv0}
\int_{\R^d}(|\log |x|| \wedge |x|^{-2}) \tau(dx) < \infty
\end{equation} and then we have
\begin{equation}\label{eq:harexpThor}
\psi_\mu(z)=i \langle z, b\rangle +\frac{1}{2}\langle \Sigma z , z \rangle  - \int_{\R^d} \left( \log \left(1  - \frac{i \langle z, x \rangle}{|x|^2} \right) +\frac{i \langle z, x \rangle}{|x|^2}\I_{\{|x| \geq 1\}}   \right) \tau(dx) ,  \quad z \in \mathbb R^d,
\end{equation}
with $b=(b_1,\dots, b_d)$ given by 
\begin{equation}\label{eq:b}
b_j=a_j- \int_{\R^d}x_j\left(\frac{1-e^{-|x|}}{|x|^2} -\frac{1}{|x|^2}\I_{\{|x| \geq 1 \}} \right)\tau(dx), \qquad j=1,\ldots, d.
\end{equation}

Conversely, let $b \in \mathbb R^d$,  $\Sigma$  be a $d \times d$ nonnegative definite symmetric matrix, and  $\tau$ a measure on $\R^d$ given by \eqref{eq:ThorinM}, where $\lambda$ is a finite measure on $S^{d-1}$ with  $\{\tau_u, \, u \in S^{d-1} \}$ a family of $\sigma$-finite measures on $[0,\infty)$ such that $\tau_u(B)$ is measurable for all $B \in \mathcal B(\R_+ )$, and further satisfying $\mathscr L(\tau_u,r)<\infty$ for all $r>0$, $u \in S^{d-1}$. Moreover, assume that the integrability condition \eqref{eq:thorconv0} holds.

Then 
 \eqref{eq:harexpThor} is the characteristic exponent of some $\mu \in T(\R^d)$ whose  L\'evy triplet is  $(a, \Sigma, \nu)$ where $a$ is the solution of \eqref{eq:b}, and $\nu$ has polar coordinates $ \left(\lambda, \{ \mathscr L(\tau_u,r)r^{-1}dr, \, u \in S^{d-1} \} \right)$.

\end{theoremEnd}

\begin{proofE}

First of all let us observe that by virtue of \eqref{eq:ThorinM} we have that \eqref{eq:thorconv0} is equivalent to
\begin{equation}\label{eq:thorconv}
\int_{S^{d-1}}\int_0^\infty (|\log s| \wedge s^{-2})\lambda(du) \tau_u(ds) < \infty.
\end{equation} 
From $\int_{S^{d-1}} \int_0^\infty \lambda(du) (r \wedge r^{-1}) k(r,u) dr <\infty$ and \eqref{eq:thorinm} we obtain through Fubini's Theorem
\begin{align}\label{eq:integrabK}
\infty &> \int_{S^{d-1}} \int_0^\infty \lambda(du) (r \wedge r^{-1}) k(r,u) dr  \nonumber  \\ &=\int_{S^{d-1}}\lambda(du)\int_0^\infty \left(\int_0^1 r e^{- r s}  dr\right) \tau_u(ds) + \int_{S^{d-1}}\lambda(du)\int_0^\infty \left(\int_1^\infty  \frac {e^{- r s}}{r} dr\right) \tau_u(ds) \nonumber \\ &= 
\int_{S^{d-1}}\lambda(du)\int_0^\infty \frac{1- e^{-s}(s+1)}{s^2} \tau_u(ds) + \int_{S^{d-1}}\lambda(du)\int_0^\infty \Gamma(0,s) \tau_u(ds)
\end{align}
with $\Gamma(\cdot, \cdot)$ the upper incomplete gamma function. The integrands in the last line of the above  are  positive, hence for the whole expression to be finite their  integrals must (separately) converge. For the first one, as $s \sim 0$ then $(1- e^{-s}(s+1))s^{-2} \sim (1-e^{-s})s^{-1} \sim 1$ and  as $s \rightarrow \infty$, it is $(1- e^{-s}(s+1))s^{-2}  \sim s^{-2}$. For the second one, recall $\Gamma(0,s) \sim -\log s, s\sim 0$ and $\Gamma(0,s) \sim e^{-s}$, $s \rightarrow \infty$. The leading convergence order of the two expressions combined is thus $|\log s|$ around 0 and $s^{-2}$ at $\infty$ and thus \eqref{eq:thorconv} and \eqref{eq:thorconv0} are proved.


Now  write the integral  \eqref{eq:LK} in polar coordinates recalling that 
$m_u$ is as in \eqref{eq:polarmease} with canonical function $k$ given by \eqref{eq:thorinm}, obtaining, again by Fubini's Theorem
\begin{equation}\label{eq:Texp1}
\psi(z)=i \langle z, a \rangle +\frac{1}{2}\langle \Sigma  z,  z  \rangle + \int_{S^{d-1}}\lambda(du)\int_{0}^{\infty}   \tau_u(ds)  \int_0^\infty \left(e^{i  r\langle z, u \rangle } -1 - i r \langle z ,  u \rangle \I_{\{ r \leq 1\}} \right)  \frac{ e^{-rs}}{r}  dr .
\end{equation}
Using the exponential and logarithmic series and interchanging summation and integration by absolute convergence, from \eqref{eq:Texp1} it follows 
\begin{align}\label{eq:Texp3}
&\int_{0}^{\infty}   \tau_u(ds)  \int_0^\infty \left(e^{i  r\langle z, u \rangle } -1 - i r \langle z, u  \rangle \I_{\{r \leq 1 \}} \right)  \frac{ e^{-rs}}{r}  dr  \nonumber  \\ & =\int_{0}^{\infty}   \tau_u(ds) \left( \int_0^\infty \sum_{k \geq 2} \frac{ (i \langle z, u \rangle)^k}{k!} r^{k-1} e^{-rs}  dr \nonumber \right.  
+ \left. i \langle z , u \rangle  \int_1^\infty e^{-r s} dr \right)  \nonumber \\ &=\int_0^\infty \tau_u(ds) \int_{0}^{\infty}   \left( \sum_{k \geq 2} \frac{ (i \langle z, u \rangle)^k}{k!}   \frac{\Gamma(k)}{s^k}  +  i \langle z , u \rangle \frac{ e^{-s} }{s} \right) 
 \nonumber \\ &= \int_{0}^{\infty}   \tau_u(ds)   \left(- \log \left(1  - \frac{i \langle z, u \rangle}{s} \right) -i \langle z, u \rangle\frac{1-e^{-s}}{s}   \right)  
\end{align}
with $\Gamma(\cdot)$ the Euler gamma function. Substituting in  \eqref{eq:Texp1} and using the the expression of $\tau$ in \eqref{eq:ThorinM} it follows 
\begin{equation}\label{eq:expThorWrongTrunc}
\psi_\mu(z)=i \langle z, a\rangle +\frac{1}{2}\langle \Sigma z , z \rangle  - \int_{\R^d} \left( \log \left(1  - \frac{i \langle z, x \rangle}{|x|^2} \right) +i \langle z, x \rangle \frac{1-e^{-|x|}}{|x|^2}   \right) \tau(dx) ,  \quad z \in \mathbb R^d.
\end{equation}
Recalling \eqref{eq:b}, from \eqref{eq:expThorWrongTrunc} we obtain \eqref{eq:harexpThor}. Notice that in view of \eqref{eq:thorconv0} the integrals on the right hand side of \eqref{eq:b} do converge. To see this one simply observes that as $|x| \rightarrow 0$ it is $ (1-e^{-|x|})|x|^{-1} \sim 1$,  and that convergence at infinity is of leading order $|x|^{-1}$.

 For the converse statement, we need only to prove the last claim, since that $\mu \in T(\R^d)$ would then follow by Theorem \ref{thm:ThorCM}. Assume then that  \eqref{eq:thorconv} holds and that
\begin{align}\label{eq:intfinass}
 \int_{S^{d-1}}\lambda(du) \int_{0}^{\infty}     \left(- \log \left(1  - \frac{i \langle z, u \rangle}{s} \right) -\frac{i \langle z, u \rangle}{s}\I_{\{s \geq 1\}}   \right) \tau_u(ds)< \infty.
 \end{align}
Using \eqref{eq:b} we can  write 
\begin{align}\label{eq:Texp3conv}
\infty &> \int_{S^{d-1}}\lambda(du) \int_{0}^{\infty}     \left(- \log \left(1  - \frac{i \langle z, u \rangle}{s} \right) -\frac{i \langle z, u \rangle}{s}\I_{\{s \geq 1\}}   \right) \tau_u(ds)+   i \langle  z, b \rangle +\frac{1}{2}\langle \Sigma z, z \rangle 
 \nonumber \\ &= \int_{S^{d-1}} \lambda(du)\int_{0}^{\infty}   \tau_u(ds)   \left(- \log \left(1  - \frac{i \langle z, u \rangle}{s} \right) -i \langle z, u \rangle\frac{1-e^{-s}}{s}   \right) + i \langle z,  a \rangle+\frac{1}{2}\langle \Sigma z, z \rangle.
\end{align}
 Substituting the inner integral of the last line of \eqref{eq:Texp3conv} in the last line of \eqref{eq:Texp3}, reading the equalities backwards, and substituting back in \eqref{eq:Texp3conv}, the resulting expression leads to \eqref{eq:Texp1}.

Hence only \eqref{eq:intfinass} remains to be shown. Observe that
based on the series for $\log(1-z)$, $z \in \mathbb C$, we have 
\begin{equation}\label{eq:convintdrift}
 \int_{1}^{\infty}     \left(- \log \left(1  - \frac{i \langle z, u \rangle}{s} \right) -\frac{i \langle z, u \rangle}{s}  \right) \tau_u(ds)=   \sum_{k \geq 2} \frac{ (i \langle z, u \rangle)^k}{k}  \int_0^\infty \tau_u(ds)s^{-k}<\infty \end{equation}
  since \eqref{eq:thorconv} implies that $c_k:= \int_0^\infty \tau_u(ds)s^{-k}<\infty$, $k \geq 2$, and furthermore $(c_k)_{k \geq 2}$ is a decreasing sequence, while the outer series is known to be convergent. 
 Also
\begin{align}\label{eq:convinitdrift2}
  \int_{0}^{1}  \left(- \log \left(1  - \frac{i \langle z, u \rangle}{s} \right) \right)\tau_u(ds) 
  = \int_{0}^{1} \log  s  \, \tau_u(ds) -  \int_{0}^{1}  \log (s-  i \langle z, u \rangle ) \tau_u(ds) .
\end{align}
The convergence of the first integral in \eqref{eq:convinitdrift2} is clear by \eqref{eq:thorconv}. Furthermore, for all $z \in \R^d$, $s \in [0,1],  u \in S^{d-1}$ 
\begin{equation}
\mbox{Re}\Big(\log (s-  i \langle z, u \rangle )\Big) =\log |s-  i \langle z, u \rangle |  \leq \frac{1}{2} \log( s^2+|z|^2) \leq  \frac{1}{2} \log( 1+|z|^2)
\end{equation}
having used the Cauchy-Schwarz inequality. Moreover, Im$ \Big(\log (s-  i \langle z, u \rangle )\Big)=\mbox{Arg}(s-  i \langle z, u \rangle)$, with Arg the principal value of the complex argument, which is bounded by $\pi$ for all $z,s,u$.   Convergence of the second integral in \eqref{eq:convinitdrift2} then follows from one last application of \eqref{eq:thorconv}. 

We therefore conclude that \eqref{eq:intfinass} holds and the proof of the theorem is finished.
 \end{proofE}

For $d=1$ the representation of Theorem \ref{thm:olddef} exactly coincides with that provided by Thorin and Bondesson. The following definition is then natural.

\begin{defn}
The measure $\tau$ in \eqref{eq:ThorinM} is the \emph{Thorin measure} associated with $\mu \in T(\R^d)$ and $(b, \Sigma, \tau)$ the \emph{Thorin triplet} of $\mu$. The radial component $\{\tau_u, u \in S^{d-1}\}$ in the polar representation of $\tau$ is its \emph{Thorin family}.  
\end{defn} 
The definition of Thorin process can be then given as  follows.
\begin{defn} A \emph{Thorin process} (or L\'evy--Thorin process) is a self-decomposable L\'evy process $X=(X_t)_{t \geq 0}$, such that $\mathcal L(X_1) \in T(\R^d)$, or equivalently if $\mathcal L(X_t) \in T(\R^d)$ for all $t > 0$. A Thorin process is uniquely determined by its Thorin triplet.
A \emph{Thorin subordinator} is an almost-surely (a.s.) increasing Thorin process on the real line, and its unit time distribution is that of a $T(\R_+)$ (GGC) distribution. 
\end{defn}
By Proposition \ref{prop:spherrep}, the Thorin measure $\tau$ is characterized by  its polar representation $(\lambda, \{\tau_u, \, u \in S^{d-1}\})$ up to multiplication by a measurable function. 

\begin{oss} In the spirit of \citet{ros:07}, it is also possible to associate to $\mu \in T(\R^d)$ a ``dual'' Thorin measure $\tau^*$  defined as the pushforward $\tau^*=p_\sharp\tau$ with $p(x)=x/|x|^2$.
Such definition would lead to the following  representation for the characteristic exponent
\begin{equation}\label{eq:harexpThordual}
\psi_\mu(z)=i \langle z, b\rangle +\frac{1}{2}\langle \Sigma z , z \rangle  - \int_{\R^d} \left( \log \left(1  - i \langle z, x \rangle \right) + i \langle z, x \rangle\I_{\{|x| \leq 1\}}   \right) \tau^*(dx) ,  \quad z \in \mathbb R^d.
\end{equation}
This is in fact the choice of \citet{lev+al:21}. Instances of this definition of Thorin measure are less common in the literature, and also not very fit to our purposes. Indeed, the slightly more appealing expression \eqref{eq:harexpThordual} compared to \eqref{eq:harexpThor} comes at the cost of  complicating the mathematics in Section 6. 
\end{oss}

\section{Path and distributional properties}\label{sec:distprop}

As a rule of thumb the order of integrability of the moments of the L\'evy measure around zero corresponds to the order of  path variation, whereas integrability at infinity to the existence of distributional moments. Similar relationships exist when considering the Thorin measure, but with the role of the extrema reversed.  This is not surprising, and can be interpreted intuitively in light of the Karamata Tauberian theorem: since the canonical function and the Thorin measure are related through a Laplace transform/inversion, the tail order of the former around 0  (respectively, at infinity) is equivalent to that of the latter at infinity (resp. around 0).  

The following results generalize \citet{buc+al:17}, Proposition 2.6,
 to the full $T(\R^d)$ class. We defer examples until later in Section \ref{sec:ThorList}.

\subsection{Path variation}

For $[a,b] \subset \R$, let $\Pi$ be the set of all the partitions  $\pi=\{a=t_0 <t_1 \leq \cdots \leq t_n=b \}$ of $[a,b]$ for $n$ varying in $\mathbb N$. The $p$-variation $V_p(f; a, b)$  of a function $f:[a,b] \rightarrow \R^d$ is defined as

\begin{equation}
V_p(f;a,b)=\sup_{\pi \in \Pi }\sum_{i=0}^n\|f(t_{i})-f(t_{i-1})\|^p. 
\end{equation}
If $X=(X_t)_{t \geq 0}$ is a stochastic process on $\R^d$, for all $[0,T] \subset  \R$ we can define the r.v. $V_p(X,  T)$ by $V_p(X,  T)(\omega):=V_p(X_t(\omega);0,T)$ 
and  one can introduce the $p$-variation index $v(X)$   of $X$ as follows
\begin{equation}
v(X)=\inf\left\{ p>0 :  V_p(X,T) <\infty \mbox{ a.s., } \forall T>0 \right\}.
\end{equation}
The Blumenthal--Getoor (BG) index $\beta(X)$ of a L\'evy process $X=(X_t)_{t\geq 0}$ with triplet $(a,0, \nu)$ on $\R^d$ is defined as
\begin{equation}\label{eq:BGLevy}
\beta(X)=\inf\left\{ p>0 :  \, \int_{\{|x| \leq 1\}}\|x \|^p \nu(dx) <\infty \right\}.
\end{equation}
 \citet{blu+get:61}, \citet{mon:72}, and \citet{bre:72} have shown that $\beta(X) \in [0,2]$ and if $X$ is without Gaussian component, its BG index coincides with the $p$-variation index, i.e. $v(X)=\beta(X)$. A process is said to be of finite variation if $\beta(X) \leq 1$, and of infinite variation otherwise. For Thorin processes, the determination of the BG index translates into a condition on the Thorin measure, along the lines of the following proposition.
\begin{theoremEnd}[proof at the end,
	no link to proof]{prop}\label{prop:fv}  Let $X=(X_t)_{t \geq 0}$ be a Thorin process, with Thorin triplet $(b, 0, \tau)$. Then its {\upshape BG} index equals
\begin{equation}\label{eq:BGthor}
\beta(X)=\inf\left\{ p>0 :  \int_{\{|x| \geq 1\}}\frac{\tau (dx)}{|x|^p} <\infty \, \right\}.
\end{equation}
In particular $X$ is of finite variation if and only if \begin{equation}\int_{\{|x| \geq 1\}}\frac{\tau (dx)}{|x|} <\infty.
\end{equation}
\end{theoremEnd}

\begin{proofE}

Let $p \in (0,2)$. Using the polar representation of the L\'evy measure of the $T(\R^d)$ distribution $\mathcal L(X_1$) and Tonelli's Theorem, for all $u \in S^{d-1}$ it holds that
\begin{align}\label{eq:gammasconv0}
&\int_{\{|x| \leq 1\}} \|x\|^p \nu(dx)=\int_{S^{d-1}}\lambda(du)
\int_0^1 r^{p-1} k(u,r)dr 
\nonumber \\ &=\int_{S^{d-1}}\lambda(du) \int_0^\infty  \left(\int_0^1 e^{-r s} r^{p-1} dr \right)  \tau_u(ds) =  \int_{S^{d-1}}\lambda(du) \int_0^\infty  \frac{\gamma(p,s)}{s^p} \tau_u(ds)  
\end{align} 
with $\gamma(p,s)$ the lower $p$-gamma incomplete function at $s>0$. Now  $\gamma(p,s) s^{-	p} \rightarrow 1/p$ when $s \rightarrow  0$ (\citet{abr+ste:48}, Section 6.5), so that  $\gamma(p,s)$ is bounded in $[0,1]$, whereas $\gamma(p,s) \uparrow \Gamma(p)$, $s \rightarrow \infty$. Therefore, for some $c_p>0$ we have
\begin{align}\label{eq:gammasconv}
&\int_{\{|x| \leq 1\}} \|x\|^p \nu(dx)< c_p  \int_{S^{d-1}}\lambda(du)\int_{(0,1)} \tau_u(ds)   + \Gamma(p) \int_{S^{d-1}}\lambda(du) \int_1^\infty s^{-p} \tau_u(ds).
\end{align}
By \eqref{eq:thorconv} the first integral converges, whereas the second does if and only if $p > \beta(X)$. 
 Taking the infimum on $p$ of the set $\{p \in (0,2)\}$ for which convergence is achieved establishes the claim.
\end{proofE}

We introduce at this point a concept connected to path variation of importance further on. For $\mu \in T(\R^d)$ with Thorin measure $\tau$, we define the measure $\tau_\circ$ as
\begin{equation}\label{eq:taucirc}
\tau_\circ(B)=\int_{S^{d-1}}\lambda(du)\tau_u(B), \qquad B \in \mathcal B(\R_+ \setminus \{0\} )
\end{equation}
which by equation \eqref{eq:thorconv0} satisfies
\begin{equation}\label{eq:thorconvcirc}
\int_0^\infty (|\log s| \wedge s^{-2}) \tau_\circ(ds) < \infty.
\end{equation}
Observe that $\tau_\circ$ is a Thorin measure with law  $\mu_\circ \in T(\R)$. 
Also, in general, when $\mu \in$  SD$(\R^d)$, we have that $k(u, 0+)<\infty$ implies $\beta(X)=0$ but the converse does not hold (a counterexample being $k(x)= -\log |x|\I_{\{|x|<e^{-1}\}}+(e x)^{-2}\I_{\{|x|\geq e^{-1}\}}, x  \in \R\setminus \{0\}$). In the former case,  $k(u, 0+)=\tau_u(\R_+)$ and furthermore in view of \eqref{eq:taucirc} it is $\tau_\circ(\R_+)=\int_{S^{d-1}} \lambda(du)k(u,0+)$.
 This justifies the following definition.

\begin{defn}\label{def:GB0} Let $X$ be an SD$(\mathbb R^d)$ L\'evy process such that $\beta(X)=0$ and \begin{equation}\label{eq:GB0} 
0<\zeta:=\int_{S^{d-1}}\lambda(du)k(u,0+) <\infty.
\end{equation}
 Then  we say that $X$ is a L\'evy process \emph{with 0 Blumenthal-Getoor index and activity $\zeta$},  and write $X \in $ BG$_0(\zeta, \R^d)$.  If $X$  is also a Thorin process we furthermore have
$\zeta=\tau_\circ(\R_+)$. 
\end{defn}

The BG$_0(\zeta, \R^d)$ class is that of the infinite activity L\'evy process closest in small time behavior to a finite activity L\'evy process, that is, a compound Poisson process. Under this analogy, the constant $\zeta$ can be understood as an intensity rate. As we will show in Section \ref{sec:sub}, BG$_0$ processes arise naturally in Thorin subordination. The simplest BG$_0(\zeta, \R)$ process is of course the $\Gamma(\zeta, \theta)$ subordinator.  As we shall see in the next sections,  BG$_0(\zeta, \R)$ processes possess a certain degree of density regularity  and enjoy particular properties in their subordination.

\subsection{Moments and  cumulants}

Let $\mu \in$ ID$(\R^d)$ such that $\hat \mu$ is multivariate complex analytical in a neighborhood of 0, and let $X \sim \mu$, $X=(X_1, \ldots , X_d)$.  Then the cumulant generating function $\kappa_X(s):=\log E[e^{ \langle s,  X \rangle}]$ can be defined (at least) for $ |s| \in (-\epsilon, \epsilon)$, $\epsilon>0$ and $\kappa_X^{n,j}=\frac{\partial^n}{\partial s_j^n } \kappa_X(0)$, $j=1, \ldots, d$, is the $n$-th cumulant of $X_j$. We have $\kappa_X(s)=\psi_X(- i s)$, $ s \in (-\epsilon, \epsilon)$, and  \begin{equation} \label{eq:cumdef}
\kappa_X^{n,j}= (-i)^n \frac{\partial^n}{\partial s_j^ n } \psi_X(0).
\end{equation} 
 Moments  and cumulants of $X$ are in a one to one mapping, with $\kappa_X^{1,j}$ and $\kappa_X^{2,j}$ coinciding with $E[X_j]$ and Var$[X_j]$ respectively. Furthermore 
\begin{equation}\label{eq:jointcumdef}
\Cov(X_i, X_j)=
 - \frac{\partial^2}{\partial s_i \partial s_j } \psi_X(0).
\end{equation} 
  In particular the $n$-th moment is finite if and only if its $n$-th cumulant is.  Therefore if $X$ has L\'evy triplet $(a, \Sigma, \nu)$, then for $j=1, \ldots, d$, it holds  that
\begin{equation}\label{eq:meanhighercum}
\kappa_X^{1,j}=a_j +\int_{\{x > 1 \}}x_j \nu(dx), \qquad 
\kappa_X^{n,j}= \Sigma_{jj} \I_{\{n=2\}} + \int_{\R^d}x_j^n \nu(dx), \quad n \geq 2,
\end{equation}
and furthermore for $i,j=1, \ldots, d$, we have
\begin{equation}\label{eq:varcovar}
\Cov(X_i,X_j)=\Sigma_{ij} + \int_{\R^d}x_i x_j \nu(dx).
\end{equation}
See \citet{sat:99}, Corollary 25.8. For Thorin distributions, the cumulants and covariances can be expressed in terms of the Thorin measure as follows.

\begin{theoremEnd}[proof at the end,
	no link to proof]{prop}\label{prop:momprop}
Let $\mu \in T(\mathbb R^d)$ be given, with Thorin triplet $(b, \Sigma,  \tau)$ and let $X \sim \mu$. 
For $n \in \mathbb N$, we have $E[|X|^n]< \infty$  if and only if 
\begin{equation}\label{eq:existsmomThor}
\int_{\{|x| < 1\}} \frac{\tau(dx)}{|x|^n}<\infty.
\end{equation}
Moreover, let  $n \in \mathbb N$ be such that  \eqref{eq:existsmomThor} holds.  Then
\begin{align}
\kappa_X^{1,j}&= b_j + \int_{\{|x| < 1\}}x_j \frac{\tau(dx) }{|x|^2}, \label{eq:meanThor}\\  \Cov(X_i,X_j)&= \Sigma_{ij}+ \int_{\R^d}x_i x_j\frac{ \tau(dx)}{|x|^4}, \qquad \mbox{ if } n = 2, \label{eq:covThor} \\ \kappa_X^{n,j}& = (n-1)!\int_{\R^d}x_j^n  \frac{\tau(dx)}{|x|^{2n}} ,  \qquad  \mbox{ if } n > 2. \label{eq:higherThor}
\end{align} 
\end{theoremEnd}

\begin{proofE}
 Let $\nu$ be the L\'evy measure of $\mu$ and $n \geq 1$. By \citet{sat:99}, Chapter 25, the existence of the absolute moments is equivalent to $\int_{\{|x| > 1\}} |x|^n\nu(dx) <\infty$. With $|x|=r$, similarly to Proposition \ref{prop:fv}, we have by Tonelli's Theorem
\begin{align}\label{eq:testMoments0}
\int_{\{|x| > 1  \}} |x|^n\nu(dx)=&\int_{S^{d-1}} \lambda(du)\int_{(1,\infty)}r^{n-1}k(u,r) dr= \int_{S^{d-1}} \lambda(du)\int_0^{\infty} \left(\int_1^\infty e^{-r s } r^{n-1}dr\right) \tau_u(ds) \nonumber \\ =&\int_{S^{d-1}} \lambda(du) \int_0^  \infty \frac{\Gamma(n,s)}{s^n}\tau_u(ds) \end{align}  with $\Gamma(n,s)$ the upper incomplete gamma function at $s$. We have that $\Gamma( n,s) s^{-n} \sim e^{-s}s^{-1}$ (\citet{abr+ste:48}, 6.5.32), such function is locally bounded, and $\Gamma(n,s) \rightarrow \Gamma(n)$ when $s \rightarrow 0$.  Thus	 there exist some $c_n , C_n>0$ such that
\begin{align}\label{eq:testMoments}
\int_{\{|x| > 1\}} |x|^n\nu(dx) < c_n \int_{S^{d-1}} \lambda(du) \int_{(0,1)} {s^{-n}}\tau_u(ds)  + C_n \int_{S^{d-1}} \lambda(du) \int_1 ^\infty {s^{-2}}\tau_u(ds).  
\end{align}
By \eqref{eq:thorconv} the second term always converges,  and the first does if and only if \eqref{eq:existsmomThor} holds true.

Now observe that for $ j,l=1, \ldots, d$, $z \in \R^d$, and $n \in \mathbb N$ it holds that
\begin{equation}\label{eq:logder}
\frac{\partial^n}{\partial z_j^n}\log \left(1-i\frac{\langle z, x \rangle }{|x|^2}\right)= -(i^n) (n-1)! \frac{ ( x_j/|x|^2)^n }{(1- i \langle z, x \rangle/|x|^2)^{n} } \end{equation}
and
\begin{equation}\label{eq:logmixder}
\frac{\partial^2}{\partial z_j \partial z_l}\log \left(1-i\frac{\langle z, x \rangle }{|x|^2}\right)=  \frac{  x_j x_l/|x|^4 }{(1- i \langle z, x \rangle/|x|^2)^{2} }. \end{equation}
 Performing the calculations in  \eqref{eq:cumdef}  
and \eqref{eq:jointcumdef} by using respectively \eqref{eq:logder} and \eqref{eq:logmixder} in \eqref{eq:harexpThor}, after having taken the derivative under integral sign, yields \eqref{eq:meanThor}--\eqref{eq:higherThor}, provided the integrals converge as assumed in \eqref{eq:existsmomThor}.
\end{proofE}
\begin{oss}\label{oss:centerrep}
If $\int_{\{|x| < 1\}}\tau(dx)  |x|^{-1}< \infty$ (the first absolute moment exists), then \eqref{eq:harexpThor} can be used with the constant function 1 replacing $\I_{\{|x|\geq 1\}}$, in which case $b$ becomes the expected value of $\mu$. See also \citet{sat:99}, p.  39.
\end{oss}

\subsection{Exponential moments and the moment generating function}\label{subsec:em}

The question of the existence of exponential moments is a classic problem in many applications of the theory of distributions. A prime example in insurance mathematics is the calculation of the Lundberg exponent, providing a bound/approximation for the ruin probability in the Cram\'er-Lundberg model (\cite{emb+al:13}, Chapter 1). Within the Thorin class, criteria sharper than those for general i.d. laws exist. 

Let $\mu$ be a probability distribution on $\R^d$. Define the critical exponent $\gamma(\mu)$ of $\mu$ as 
\begin{equation}
\gamma(\mu)=\sup\left\{p \geq 0 : \, \int_{\R^d} e^{p |x|} \mu(dx) <\infty   \right\} \geq 0.
\end{equation}
Then $\hat \mu$ is analytical in a neighborhood of zero -- and a moment generating function (m.g.f.) exists -- if and only if $\gamma(\mu) >0$. The key property for Thorin processes is that the infimum of the support of $\tau_\circ$ coincides with the critical exponent of $\mu$.

\begin{theoremEnd}[proof at the end,
	no link to proof]{prop}\label{prop:existsmomThor} Let $\mu \in T(\R^d)$. 	 
 Then \begin{equation}
\gamma(\mu)=\inf_{[0, \infty)}  \mbox{{\upshape supp}}\, \tau_\circ \geq 0
\end{equation} with $\tau_\circ$ 
 as in \eqref{eq:taucirc}. 
 Furthermore, if $\gamma(\mu)>0$ then $\int_0^\infty e^{\gamma(\mu) |x|}\mu(dx)<\infty$ if and only if $\tau_\circ$ does not have an atom in $\gamma(\mu)$.
\end{theoremEnd}

\begin{proofE}
We can ignore the Gaussian component of $\mu$, since all of its exponential moments exist. Let $\nu$ be the L\'evy measure of $\mu$ and set $\gamma:=\gamma(\mu)$. Let us briefly recall that $\mathscr L( \tau_\circ^\gamma, r) \rightarrow 0$ as $r \rightarrow \infty$ if and only if  $\tau_\circ^\gamma(\{0\})=0$. 
Integrating by parts and by means of the substitution $u=r s$ one obtains
\begin{equation}
\mathscr L( \tau_\circ^\gamma, r)=r \int_0^\infty e^{-r s}F_{\tau_\circ^\gamma}(s) ds=\int_0^\infty e^{-u}F_{\tau_\circ^\gamma}(u/r) du, \qquad r>0.
\end{equation}
 Taking by dominated convergence  the limit as $r \rightarrow \infty$ inside the integral, we conclude that the left hand side tends to 0 if and only if $0=F_{\tau_\circ^\gamma}(0+)=\tau_\circ^\gamma(\{0\})=\tau_\circ(\{\gamma\})$.

 From  \citet{sat:99}, Corollary 25.8, existence of the exponential moments above is equivalent to $\int_{\{|x| >1 1\}}^\infty e^{|x|\theta}\nu(dx) < \infty$. 
 Let $\theta  \geq 0 $. Applying Tonelli's Theorem and the convolution rule of the Laplace transform, 
it holds that \begin{align}\label{eq:expintconv}\int_{\{|x| > 1\}}^\infty e^{|x|\theta}\nu(dx)  &= \int_{S^{d-1}}\lambda(du)\int_1^\infty   \frac{e^{\theta r}}{r} \mathscr L(\tau_u,r)  dr =  \int_1^\infty   \frac{e^{\theta r}}{r} \mathscr L(\tau_\circ,r) dr \nonumber \\ &= \int_1^\infty   \frac{e^{\theta r}}{r} \mathscr L(\tau^\gamma_\circ * \delta_\gamma,r) dr = \int_1^\infty   \frac{e^{(\theta-\gamma) r}}{r} \mathscr L(\tau^\gamma_\circ, r) dr.
\end{align}
If $\theta<\gamma$ this integral always converges because the Laplace transform is a decreasing function. If $\gamma=\theta$  the integrand is  $o(r^{-1})$ for large $r$ if and only if $\mathscr L( \tau_\circ^\gamma, r) \rightarrow 0$ when $r \rightarrow \infty$ which, in combination to what has been argued before, proves the last claim.
 
 It remains to be ascertained divergence for $\theta>\gamma$. But this is clear since for all $ a>0$ we have
\begin{equation}
\mathscr L( \tau_\circ^\gamma, r) \geq \int_0^a e^{- r s}\tau_\circ^\gamma(ds) \geq e^{-r a} \tau_\circ^\gamma([0,a])> e ^{- r a}  c_a
\end{equation} 
for some $c_a>0$, since $0 \in $ supp$\, \tau^\gamma_\circ $. Therefore $\infty \geq \liminf_{r \rightarrow \infty} \mathscr L(\tau^\gamma_\circ, r)e^{a r} >0$ and the last integral in \eqref{eq:expintconv} must diverge.
 \end{proofE}


The main insight of Proposition \ref{prop:existsmomThor} is that $\hat \mu$ can be studied through some of the measure-theoretic properties of $\tau$.  The univariate case is of obvious importance. The following corollary is a consequence of Proposition \ref{prop:existsmomThor}  and the standard theory of c.f.s, see \citet{luk:70}.
\begin{theoremEnd}[proof at the end,
	no link to proof]{cor}\label{cor:analchar} Let $\mu \in T(\R)$ and $\gamma_+=\inf_{[0, \infty]} \mbox{{\upshape supp}} \,  \tau_+$, $\gamma_-=\inf_{[0, \infty]}  \mbox{{\upshape supp}} \,  \tau_{-}$. Then $\mu$ is analytical 
if and only if $\gamma_+, \gamma_- >0$,  in which case $\hat \mu$ has convergence strip  $\mathcal S=\{ z \in \mathbb C : \mbox{ \upshape{Im}}(z) \in (-\gamma_+, \gamma_-) \}$. Furthermore,  $\hat \mu$ extends analytically on the corresponding boundary $\partial^\pm \mathcal S=\{z \in \mathbb C :   \mbox{ \upshape{Im}}(z)=\pm \gamma_\mp \}$ of $\mathcal S$ if and only if $\tau_{\pm}(\{\gamma_\pm\})  = 0$ . In any case, $\hat \mu$ can be analytically continued to an entire function if and only if $\mu$ is a Gaussian distribution.
\end{theoremEnd}

\begin{proofE}
Convergence of the Fourier transform $\hat \mu$ in the strip $\mathcal S$ is the same as the convergence of the bilateral Laplace transform in the vertical strip with abscissae $-\gamma_-$ and $\gamma_+$ which is in turn equivalent to the existence of the m.g.f. in $(-\gamma_-, \gamma_+)$. Hence the first assertion follows from Proposition \ref{prop:existsmomThor}. Also by Proposition 4 such convergence can be extended on  $\partial^+ \mathcal S$ (respectively,  $\partial^- \mathcal S$) if $\tau_+$ (respectively,  $\tau_-$) does not have an atom in $\gamma_+$ (respectively, $\gamma_-$). The last statement is clear.
\end{proofE}

We shall shortly see that finiteness of the critical exponential moment plays a primary role in determining the probability density tails of Thorin distributions.

\subsection{Univariate densities: properties and asymptotics}\label{subsec:ud}

There is a known connection between the analytical properties of the canonical function of a self-decomposable distribution and the regularity and unimodality properties of the corresponding i.d. distribution. As expected, when taking into account  Thorin distributions, equivalent conditions can be formulated in terms of the Thorin measure.

\begin{theoremEnd}[proof at the end,
	no link to proof]{prop}\label{prop:Thorfinite}
Let $\mu \in T(\R)$ with Thorin triplet $(b, \sigma, \tau)$. Setting $b_0:=\int_{\{|x|  \geq 1\}}x^{-1} \tau(dx) \leq \infty$ denote by $f$ the p.d.f. of $\mu$ and let $\zeta_\pm=\int_0^\infty \tau_\pm(dx)$, $\zeta_\pm \in \R_+ \cup \{\infty \}$ and $\zeta:=\zeta_+ + \zeta_-$ with $\infty+\infty:=\infty$. Notice that $\zeta<\infty$ corresponds to $\mu \in$ {\upshape BG}$_0(\zeta, \R)$.   The following hold true:
\begin{itemize}
\item[(i)] if either $\zeta=\infty$ (i.e. $\mu \notin $ {\upshape BG}$_0(\zeta, \R)$ for all $\zeta > 0$) or $\sigma>0$ then $f \in C^{\infty}(\R)$;
\item[(ii)] if $\sigma=0$ and $1 <\zeta < \infty$ then $f\in C^{n-1}(\R)$, and $ f \notin C^{n}(\R)$ with $n=\left \lfloor{\zeta }\right \rfloor  $;  
\item[(iii)] if $\sigma=0$ and $\zeta \leq 1$ then $f$ is continuous, except, if $b_0 <\infty$, possibly at $b_*=b-b_0$. However,  $h(x):=(x-b_*)f(x)$, $x \neq b_*$, $h(b_*)=0$, satisfies $h \in C(\R)$.  
\item[(iv)] $\mu$ is strictly unimodal i.e. there is $m \in \R$ such that $f$ is increasing on $(-\infty, m)$ and decreasing on $(m, \infty)$. 
\end{itemize}
\end{theoremEnd}

\begin{proofE} 
Theorems, equations and the classification henceforth mentioned without reference are those in \citet{sat+yam:78}. Points $(i)-(ii)$ follow from Theorem 1.2 (due to Wolfe and Zolotarev) and the Bernstein representation \eqref{eq:thorinm} of $k$. 
Recalling the relationship \eqref{eq:b} between the L\'evy linear characteristic $a$ and the Thorin  one $b$,  under the assumption $b_0<\infty$ it holds that
\begin{align}
a-\int_{-1}^{1}x \frac{k(x)}{|x|}dx&=a -\int_0^{1}dr \left(\int_0^\infty e^{-s r}\tau_+(ds)-\int_0^\infty e^{-s r}\tau_-(ds) \right) \nonumber \\ &=a - \int_{-\infty}^{\infty}\frac{1-e^{-|x|}}{x}\tau(dx) = b - b_0=b_*
\end{align}
with an application of Fubini's Theorem. 
The first summand of the above is the drift $\gamma_0$ in Equation (1.6) (upon replacing the indicator truncation function to the one in \cite{sat+yam:78}), and we can apply the theorems in \citet{sat+yam:78} with $\gamma_0=b_*$. 
In particular now  $(iii)$ also follows from Theorem 1.2. 
Moreover, neither $\mu$ nor $\check{\mu}$ is of type I$_4$ because $k_\pm$ is c.m. and not constant, hence strictly decreasing, so that if $k(0+)=\lambda$ then $k(x)<\lambda$ for all $x>0$. Now $(iv)$ is implied by  Theorem 1.4. 
\end{proofE}

Using the results in \citet{sat+yam:78}, several other properties on the mode, location and regularity on the left and right sides of $b_*$, follow from the analysis of $\zeta_\pm$. We omit these for the sake of brevity as they follow easily from considerations analogue to those in Proposition \ref{prop:Thorfinite}.

\bigskip
We turn now to the asymptotic properties of univariate densities. Some definitions are needed.
For $\gamma \geq 0$, a distribution $\mu$ on $\mathbb R$  belongs to the class $L(\gamma)$   if for all $x,y \in \mathbb R$ we have $\overline \mu(x)>0$ and 
 \begin{equation}\label{eq:convlclass}\frac{\overline {\mu}(x+y)}{ \overline \mu(x)}=e^{- y \gamma}, \qquad x \rightarrow + \infty .\end{equation} The distribution $\mu$ is instead said to be  \emph{convolution equivalent}, or of class $\mathcal S(\gamma)$, if $\hat \mu(-i \gamma)<\infty$, $\mu \in L(\gamma)$, and
\begin{equation}\label{eq:yammain}
\lim_{x\rightarrow \infty}\frac{\overline{\mu*\mu}(x)}{\overline{\mu}(x)}= 2 \hat \mu (  -i\gamma)<\infty.
\end{equation}
See  e.g. \citet{wat:08}, Equations (1.1) and (1.2). For  distributions $\mu$ with support unbounded from below convolution equivalence at $x \rightarrow -\infty$ is defined as the convolution equivalence of $\check \mu$.

The connection between convolution equivalence and the tail properties of i.d. probability distributions is exposed by a  deep theorem (\citet{wat:08}, \citet{pak:04}, \citet{sgi:90}) with a rather interesting history (see \citet{wat:08}), asserting, among the other things,  that  under convolution equivalence of order $\gamma$ the complementary distribution functions of i.d. laws are asymptotically equivalent to their L\'evy tail functions, with proportionality constant the m.g.f. at $\gamma$.


Under a 
 regular variation assumption, convolution equivalence in the Thorin class is equivalent to the existence of the critical exponential moment.

\begin{theoremEnd}[proof at the end,
	no link to proof]{thm}\label{thm:subexp}
Let $\mu \in T(\R)$, with Thorin measure $\tau$ and  $ 0 \leq \gamma=\inf_{[0,\infty)} \mbox{{\upshape supp}} \, \tau_+$. If  $\hat \mu(-i \gamma)<\infty$ and $F_{\tau^\gamma}$ is regularly varying in 0 then   $\mu \in \mathcal S(\gamma)$. 
\end{theoremEnd}

\begin{proofE}
Define the L\'evy tail distribution $\nu_1$ as the probability law with density given by
 \begin{equation}
\nu_1(x):=\frac{ \nu(x)  }{\overline \nu(1)}\I_{\{x>1\}}  =\frac{ k(x)  }{ x \overline \nu(1)}\I_{\{x>1\}}. \end{equation} 
  By \citet{wat:08}, Theorem B,  $\mu \in \mathcal S(\gamma)$ if and only if  $\nu_1 \in \mathcal S(\gamma)$ which we will prove.

Consider the positive part $k_+$ of the canonical function $k$ of $\mu$, which we denote for notational convenience and with abuse of notation, still $k$. Since $k$ is c.m., by a slight extension\footnote{In the Proposition mentioned, the assumption of the Laplace transform being that of a probability measure can be lifted upon restricting the abscissa of convergence to some $s_0 >\gamma$. The proof carries through by appropriately modifying the bound in (2.3) using the Laplace transform computed at $s_0$ instead of the probability tail.} of  \citet{tor:24}, Proposition 1, $(i)$,  we have that \eqref{eq:convlclass} is satisfied by $k$. Thus  
\begin{align}\label{eq:longtailedTGS}
& \lim_{x \rightarrow \infty} \frac{\overline{ \nu_1}(x+y)}{\overline{ \nu_1}(x)}= \lim_{x\rightarrow \infty}\left( \frac{x}{x+y} \right)  \frac{k(x+y)} {k(x)}=e^{-\gamma y}, \qquad y \in \R
\end{align}
which shows  $ \nu_1 \in  L(\gamma)$. 
With the same formal computations in \citet{tor:24}, Theorem 1, one obtains
\begin{align}\label{tailslim}
\lim_{x\rightarrow \infty}\frac{\overline{\nu_1* \nu_1}(x)}{\overline{\nu_1}(x)}&=  \frac{1}{\overline  \nu(1)} \int_1^{x/2}\left( \frac{x} {(x-z)z}\right)\frac{k(x-z)k(z)}{k(x)}  dz.
\end{align}
Now for  $z \in (1,x/2)$ we have $x/(x-z) \leq 2$. Moreover, as a consequence of \citet{tor:24} Proposition 1, $(ii)$, for all fixed $ z \geq 1$ the ratio $k(x-z)/k(x)$ is decreasing in $x$ for all $x >z$  so that $k(x-z)/k(x)\I_{\{1 \leq z \leq x/2\}}$ has a maximum in $x=2 z$. Hence
\begin{align}\label{eq:int}
\left( \frac{x} {(x-z)z}\right) &\frac{k(x-z)k(z)}{k(x)}\I_{\{1 \leq z \leq x/2 \} }  < 
\left(\frac{2} {z}\right)\frac{k(z)^2}{k(2z)} \I_{\{ z \geq 1  \} }=
\left(\frac{2} {z}\right) \frac{(k(z) e^{z \gamma})^2}{k(2z)e^{2 z \gamma}}   \I_{\{ z \geq 1  \} }
\nonumber \\  
&=\left(\frac{2} {z}\right) \frac{(k^\gamma(z))^2} {k^\gamma(2 z)}   \I_{\{ z \geq 1  \} }
  \end{align}
  with $k^\gamma(z):=k(z) e^{z \gamma}=\mathscr L(\tau^\gamma_+, z)$, and $\inf_{[0,\infty)}$ supp$\,\tau^\gamma_+=0$ by convolution. 
 But now  being $F_{\tau^\gamma}$ regularly varying in $0$ it is implied by Karamata's Tauberian  Theorem in the version of \citet{fel:71}, theorems XIII.5.2 and XIII.5.3, that $k^{\gamma}$ is regularly varying at infinity, whence  $k^{\gamma}(z)/k^{\gamma}(2z) \rightarrow 2^{\rho}$, $\rho<0 $, $z \rightarrow \infty$.
 We conclude that \begin{equation}
\frac{(k^\gamma(z))^2}{k^\gamma(2 z)} \rightarrow 0, \qquad z \rightarrow \infty
\end{equation} 
and thus the right hand side of \eqref{eq:int} is integrable.  
We can then apply the Dominated Convergence Theorem and continue \eqref{tailslim}  by taking the limit inside the integral.
Using  $k(x-z) \sim e^{\gamma z}k(x)$ when $x\rightarrow \infty$, we find
\begin{align}\label{eq:final}
\lim_{x \rightarrow \infty}&\frac{\overline{\nu_1*\nu_1}(x)}{\overline{\nu_1}(x)}=\frac{2  }{\overline \nu(1)}  \int_1^{\infty}\lim_{x\rightarrow \infty}  \left( \frac{x}{(x-z)z}\right)\frac{k(x-z) k( z) }{k( x)}\I_{\{1 \leq z \leq x/ 2 \} } dz \nonumber \\&= \frac{2 }{\overline \nu(1)} \int_1^\infty \frac{e^{\gamma z} k(z) }{z} dz= 2 \hat {{\nu_1}}( - i \gamma) 
\end{align}
which is moreover finite because by assumption so is $\hat \mu(- i \gamma)$. This completes the proof that $ \nu_1 \in \mathcal S(\gamma)$, and hence that $\mu \in \mathcal S(\gamma)$. 
\end{proofE}

Combining  Theorem \ref{thm:subexp} with Proposition \ref{prop:existsmomThor} we conclude that a test for convolution  equivalence is the presence of an atom in the infimum of the supports of the Thorin family.

\begin{cor}\label{cor:testCE} Let $\mu \in T(\R)$ be a Thorin distribution as in Theorem \ref{thm:subexp}, with  Thorin family $\{\tau_+,\tau_-\}$, and let $\gamma_\pm>0$ be the critical exponents of $\mu_\pm$.  We have $\mu_\pm \in \mathcal S(\gamma_\pm)$	
if and only if $\tau_\pm(\{\gamma_\pm\})=0$
\end{cor}
As another corollary we establish the tail  equivalence of p.d.f.s and L\'evy densities of Thorin laws.
\begin{theoremEnd}[proof at the end,
	no link to proof]{cor}\label{cor:tails}
Let $\mu \in T(\R)$ satisfy the assumptions of Theorem \ref{thm:subexp} and let $\gamma_\pm=\inf_{[0, \infty)} \mbox{\upshape supp} \, \tau_\pm $ 
 The p.d.f. $f$ of $\mu$ satisfies
\begin{equation}\label{eq:pasympt}
f(x) \sim  \widehat {\mu_\pm}(- i \gamma_\pm ) \frac{k(x)}{|x|}, \qquad x \rightarrow \pm \infty.
\end{equation}
\end{theoremEnd}
\begin{proofE}
We only treat the case $x \rightarrow \infty$, since for $-x$ large in absolute value the claim would follow by considering $\check \mu$. Again we drop the subscript $+$ for clarity. By \citet{bin+al:89}, Theorem 1.4.1, \emph{(iii)}, 
we have $k^\gamma(x)=x^{\rho}\ell(x)$, $\rho \in \R$, and  $\ell$ is a slowly varying function at infinity. Furthermore,  by \citet{fel:71}, VIII.9, Theorem 1, (a) with $\nu^\gamma(x)=k^\gamma(x)/x$, we also have
\begin{equation}\label{eq:nuRV}
  \overline {\nu^\gamma}(x) \sim -  \frac{k^{\gamma} (x)}{\rho} \sim -\frac{ x^{\rho}\ell(x)}{\rho},  \qquad x \rightarrow +\infty.
\end{equation}
Now, $\mu^\gamma \in \mathcal S(0)$ by Theorem \ref{thm:subexp}; in this case \citet{wat:08}, Theorem B, (3), states that $\overline {\mu^\gamma}(x) \sim \overline {\nu^\gamma}(x)$, so that based on \eqref{eq:nuRV} it also follows 
$\overline{\mu^\gamma}(x) \sim  -x^{\rho}\ell(x)/\rho$. We observe now that by Proposition \ref{prop:Thorfinite}, \emph{(iv)}, $f$ is ultimately monotone decreasing, so that we can apply the Monotone Density Theorem in its version for the tail function (\citet{bin+al:89},  theorems 1.7.2, 1.7.2b and the related discussion), stating that if $U(x)=\int_x^\infty u(y)dy$ with $u$ ultimately monotone, and $U(x)\sim  c x^\rho \ell(x)$, $x \rightarrow \infty$, $\ell(x)$ slowly varying, $\rho <0$, then $u(x) \sim - c \rho x^{\rho-1} \ell(x)$, $x \rightarrow \infty$. In the present  case from \eqref{eq:nuRV} we set $c=-1/\rho$  we have \begin{equation}\label{eq:pgammaasy}
f^\gamma(x) \sim   x^{\rho-1} \ell(x) \sim x^{-1}k^{\gamma}(x), \qquad x \rightarrow \infty.
\end{equation} with  $f^\gamma$ the probability density of $\mu^\gamma$. Finally, $f^\gamma$ is also equal to the exponential tilting of $f$, i.e. \begin{equation}\label{eq:tilt}
f^\gamma(x)=\frac{e^{\gamma x}f(x)}{\hat \mu(-i \gamma)}
\end{equation}
for all $x$.
Expressing $f$ for $f^\gamma$ in \eqref{eq:tilt}, using \eqref{eq:pgammaasy} and recalling $k^\gamma(x)=e^{\gamma x}k(x)$, \eqref{eq:pasympt}  readily follows.
\end{proofE}

 Regular variation of $F_{\tau^\gamma}$ in 0 -- or equivalently by Karamata's Theorem that of $k^\gamma$ at $\infty$ -- is a rather minimal 
assumption commonly satisfied by popular Thorin distributions such as those analyzed in the next section, where we discuss examples and counterexamples of Theorem \ref{thm:subexp} and of its corollaries.

\section{Infinitely divisible distributions and L\'evy processes in their Thorin representation}\label{sec:ThorList}

We analyze some known and less known self-decomposable L\'evy distributions and models in their Thorin structure through the lenses of the theory developed in sections \ref{sec:Thordists} and \ref{sec:distprop}. 

\subsection{Multivariate gamma distributions}\label{subsec:gamma}
For a continuous function $\theta:S^{d-1} \rightarrow \R_+$ 
bounded away from zero, consider the canonical function $k$ and corresponding Thorin family $\{\tau_u , \, u \in S^{d-1} \}$  given by
\begin{equation}\label{eq:multigam}
k(u,r)= e^{-r \theta(u)}, \qquad \tau_u=\delta_{\theta(u)}.
\end{equation}
This law is the multivariate gamma distribution studied in \citet{per+al:12}. In particular for $ x \in \R^d,$ $x \neq 0$, $a,b>0$  letting $\theta(u)=b/|x|>0$ for all $u \in S^{d-1}$, and $\lambda(du)=a \delta_{x/|x|}(du)$ a simple variable change shows that $\mu \in \Gamma_e(\R^d)$ with  $U \sim \Gamma(a, b)$ (see \citet{bns+al:06}, page 27). We notice that being $\tau$ purely atomic, based on Proposition \ref{prop:existsmomThor}, Theorem \ref{thm:subexp} and Corollary \ref{cor:testCE}, multivariate gamma variables do not admit an exponential $\gamma(\mu)$ moment, nor are they convolution equivalent.  All these facts are well-known or easily checked in the case $d=1$, when accordingly \eqref{eq:pasympt} is not satisfied by the usual univariate gamma p.d.f.

If $X=(X_t)_{t \geq 0}$ is a L\'evy process such that $\mathcal L(X_1)$ is a multivariate gamma law of the form \eqref{eq:multigam}, we  clearly have $\int_{S^{d-1}} \lambda(du)\int_1^\infty s^{-p} \tau_u(ds) < \infty$ for all  $p>0$ so that by Proposition \ref{prop:fv} it follows $\beta(X)=0$ as already noted for gamma subordinators. Moreover the computation in  \eqref{eq:GB0} indicates that $X \in$ BG$_0(\lambda(\R^d), \R^d)$. In the case $d=1$, denoting $\lambda(\pm 1)=\lambda_\pm$, $\theta(\pm 1)=\theta_\pm$, $\lambda_\pm, \theta_\pm >0$, $X$ reduces to the familiar bilateral gamma process Bil$\Gamma(\lambda_+, \lambda_-, \theta_+, \theta_-)$ and when  $\lambda_{-}=0$ we obtain a gamma $\Gamma( \lambda_+, \theta_+)$ subordinator. From Proposition \ref{prop:momprop} we find the usual L\'evy cumulant expressions $\kappa^{n}_{X_t}=t (n-1)!(\lambda_+ \theta_+^{-n} +(-1)^n\lambda_-\theta_-^{-n})$, $n \geq 1$ (e.g. \citet{kuc+tap:08}, Equation (2.8)).  Finally in the univariate case $ \tau_\circ(\R)=\lambda_+ +\lambda_-$, and Proposition \ref{prop:Thorfinite} recovers the regularity and boundedness properties of the (bilateral) gamma process p.d.f.

\subsection{Stable laws}\label{subsec:stable}

Following \citet{sat:99}, Theorem 14.3, the $\alpha$-stable distribution class S$_\alpha(\R^d)$,  $\alpha \in (0,2)$,  can be characterized as the subclass of SD$(\R^d)$ of laws  without Gaussian component and canonical function  given by $k(u,r)=r^\alpha$.
Then, since 
 \begin{equation}\label{eq:invfracLapl}
r^{-\alpha} =\frac{1}{\Gamma(\alpha)} \int_0^\infty e^{-  r s} s^{\alpha-1} ds , \qquad r>0
\end{equation}
 the Thorin family of $\mu \in $ S$_\alpha(\R^d)$ is given by $\{{\tau_u,  \, u \in S^{d-1}}\}=\{\tau \}$ with 
\begin{equation}\label{eq:tmstable}
  \tau(ds)= \frac{s^{\alpha-1}}{\Gamma(\alpha)} ds.
\end{equation}
 The well-known stable distributional properties are immediate taking the Thorin view. Since $0 \in$ supp$\, \tau$ it follows from Proposition \ref{prop:existsmomThor} that $\hat \mu$ is not analytical in a neighborhood of 0, and thus $\mu$ does not admit a m.g.f. We have $\int_{S^{d-1}}\lambda(du)\int_0^1 s^{-1}\tau(ds)<\infty$ if and only if $\alpha \in (1,2)$ which is  the well-known condition of existence of the expectation, whereas $\int_{S^{d-1}}\lambda(du) \int_0^1 s^{-n} \tau(ds)=\infty$ for all $n \geq 2, \alpha \in (0,2)$. Also
\begin{equation}
\int_{S^{d-1}}\lambda(du) \int_{1}^\infty \frac{\tau(ds)}{s^p} <\infty
\end{equation}
for all $p > \alpha$ so by Proposition \ref{prop:fv} if $X=(X_t)_{t \geq 0}$ is a stable process we have $\beta(X)=\alpha$, and $X$ is of finite variation if and only if $\alpha \leq 1$. Theorem \ref{thm:subexp} also applies with $\gamma=0$ (long-tailedness of stable distributions), and Corollary \ref{cor:tails} establishes the asymptotic equivalence between the L\'evy and probability density tails with proportionality constant 1.

\subsection{Tempered stable laws}\label{subsec:TS}

Exponentially tempered stable distributions were introduced by \citet{kop:95} refining an original idea of \citet{man+sta:95} of truncating the tails of a L\'evy measure to ensure existence of higher moments. 
A general notion of completely monotone tempering was later proposed by \citet{ros:07} along the following lines. We define $\mu \in$ TS$_\alpha(\R^d)$, $\alpha \in (0,2)$  if $\mu \in$ SD$(\R^d)$ is without Gaussian part  and has canonical function
\begin{equation}\label{eq:kts}
k(u,r)=\frac{q(u,r)}{r^\alpha}
\end{equation}
where for all $u \in S^{d-1}$, $q: \R_+ \times S^{d-1} \rightarrow \R_+$ is a c.m. function of $r$ such that $q(r,u) \rightarrow 0$ as $r \rightarrow \infty$. 

 It is clear that TS$_\alpha(\R^d) \subset T(\R^d)$ since $q(u,r)r^{-\alpha}$ is c.m., it being product of c.m. functions.  What we show in the following is that  the corresponding Thorin family has an interpretation within the theory  of fractional calculus. 

\smallskip

Let  $a>0$ and denote $L^1(D)$  the Banach space of integrable functions on the measure space $(D, \mathcal B(D), \mbox{Leb}(D))$, $D \subseteq [0, \infty)$. The fractional integral of order $a>0$ and base point  $b\geq 0$ is the  linear operator $I_{b}^{a}:L^1(D) \rightarrow  L^1(D)$ given by
\begin{equation}\label{eq:fracint}
 I_{b}^{a}f(x)=\frac{1}{\Gamma(a)}\int_b^x (x-t)^{a-1}f(t) dt, \qquad f \in L^1(D), \quad x>b, \quad x,b \in D.
\end{equation}
Classic references on fractional calculus are \citet{kil+al:93} and \citet{gor+al:14}.
As in \citet{car+com:20}, Remark 2 it is possible to extend this definition to the Banach space $\mathcal M(D)$ of the Radon measures on $(D, \mathcal B(D))$ by letting $I_{b_+}^{a}: \mathcal M(D) \rightarrow L^1(D)$ to be 
\begin{equation}\label{eq:fracintmeas}
 I_{b}^{a}\rho(x)=\frac{1}{\Gamma(a)}\int_b^x (x-t)^{a-1}\rho(dt), \qquad \rho \in \mathcal M(D), \quad x>b, \quad x,b \in D.
\end{equation}
For $a=1$, \eqref{eq:fracint} and \eqref{eq:fracintmeas} reduce respectively to the Riemann-Stieltjes integral of a function and the Lebesgue integral of a Radon measure, and in this sense $  I_{b}^{a}$ extends the usual integral operator. 
The Thorin family of a TS$_\alpha(\R^d)$ law is constituted by the $\alpha$-fractional integrals of the measures in the Bernstein representation of the tempering function.

\begin{theoremEnd}[proof at the end,
	no link to proof]{prop}\label{prop:thoriintmeas} Let $\mu \in $ {\upshape TS}$_\alpha(\mathbb R^d)$, 
and let  $\{ Q_u , \; u \in S^{d-1} \}$ be the family of measures such that
\begin{equation}
 \label{eq:bernrep}q(r,u)=\int_0^\infty e^{-s \, r} Q_u(ds).
\end{equation} 
Then $\mu$  has Thorin measure $\tau$ having polar representation $(\lambda, \{\tau_u, \, u \in S^{d-1}\})$ with   
\begin{equation} \label{eq:tauTS}
\tau_u(ds)=I_{0}^{\alpha}Q_u(s) ds.
\end{equation}
\end{theoremEnd}
\begin{proofE}
Based on \eqref{eq:tmstable} it is possible to reinterpret  \eqref{eq:fracintmeas} as \begin{equation}\label{eq:intconv}
 I_{b}^{a}\rho=\rho * \tau_\alpha.
\end{equation}  
with $\tau_\alpha$ the Thorin measure of some $\mu \in$ S$_\alpha(\R^d)$.
We then have, using  the convolution rule of the Laplace transform  and \eqref{eq:intconv}, that
 \begin{align}\label{eq:LM2}
k(u,r)= \frac{q(r,u)}{r^\alpha}&= \int_0^\infty e^{- s r } ( Q_u * \tau_\alpha)(s) ds = \int_0^\infty e^{- s r } I_{0}^{\alpha} Q_u(s) ds.
\end{align}
and the statement follows by injectivity of the Laplace transform. 
\end{proofE}

We compute the Thorin measures in two TS$_\alpha(\R^d)$ distribution subclasses.

\begin{esmp}\label{ex:thoexptemp}\emph{Classical tempered stable}
CTS$_\alpha(\R^d)$ \emph{laws}.
The CTS$_\alpha(\R^d) \subset$ TS$_\alpha(\R^d)$ class  of exponentially tempered stable laws is defined by letting in the TS$_\alpha(\R^d)$ L\'evy representation 
\begin{equation}
q(u,r)=e^{-r \theta(u)}
\end{equation}
with $\theta: S^{d-1}\rightarrow \R_+$ continuous or bounded. Since $Q_u=\delta_{\theta(u)}$,
performing the fractional integration in \eqref{eq:tauTS} yields  
\begin{equation}\label{eq:CTStm}
\tau_u(ds)=  \frac{(s - \theta(u))^{\alpha-1}}{\Gamma(\alpha)}\I_{\{s>\theta(u) \}} ds.
\end{equation}
The corresponding characteristic exponents in the case $d=1$ are well-known. We verify that \eqref{eq:CTStm} is correct in the TS$_\alpha(\R_+)$ case, i.e. we choose $c(1)=c>0$, $\lambda(du)=\lambda \delta_1(du)$, $\lambda>0, \theta(1)=\theta>0$. As customary, we neglect the indicator function in \eqref{eq:harexpThor}, which by finiteness of the expectation  is possible upon a modification of the constant $b$ in the Thorin--Bondesson representation.  For $\alpha \neq  1$, we obtain from Fubini's Theorem and uniform convergence of the binomial series
\begin{align}\label{eq:thorinsub}
\psi_\mu(z)&=-\lambda \int_0^{\infty}\left(\log \left(1 -\frac{iz}{s}\right)+ \frac{iz}{s}\right)\tau(d s) = \frac{\lambda}{\Gamma(\alpha)} \int_{\theta}^\infty \sum_{k=2}^\infty\frac{(iz/s )^k}{ k}(s-\theta)^{\alpha-1}ds \nonumber \\ &= \frac{\lambda}{\Gamma(\alpha)} \sum_{k=2}^\infty\frac{(iz)^k}{k}\int_0^\infty x^{\alpha-1}(x+\theta)^{-k}dx  \nonumber  \\&= \frac{\lambda}{\Gamma(\alpha)}  \sum_{k=2}^\infty\frac{(iz)^k}{k}\frac{1}{\Gamma(k)}\int_0^\infty e^{-\theta y} y^{k-1}\left(\int_0^{\infty}  e^{-x y} x^{\alpha-1} dx \right) dy =  \lambda \theta^{\alpha}\sum_{k=2}^\infty\frac{(iz/ \theta )^k}{ k} \frac{\Gamma(k-\alpha)}{\Gamma(k)} \nonumber  \\& =\lambda \Gamma(-\alpha)\left((\theta-i z)^{\alpha} - iz \alpha \theta^{\alpha-1} - \theta^{\alpha} \right) 
\end{align}
which is the well-known expression for the characteristic exponent CTS$_\alpha(\R_+)$ distributions, e.g. \citet{con+tan:03}, Proposition 4.2. 
 If $\alpha=1$ we recognize instead in the third line of \eqref{eq:thorinsub} a logarithmic series, so that
\begin{align}\label{eq:thorinsub2}
&-\lambda \int_0^{\infty}\left(\log \left(1 - \frac{iz}{s}\right)+ \frac{iz}{s}\right)\tau(ds) = \lambda \theta \sum_{k=2}^\infty\frac{(iz/ \theta )^k}{ k(k-1)} = \lambda \left(  \left( \theta-i z\right) \log\left(1-\frac{i z}{\theta} \right) + iz \right) 
\end{align} which is again the right expression.

As in the stable case, the commonly known fact that CTS$_\alpha(\R^d)$ distributions are of finite variation for $\alpha<1$ and infinite variation for $\alpha \in [1,2]$ is readily recovered  from Proposition \ref{prop:fv}. Regarding the moments, CTS$_\alpha$ laws are  analytical distributions whose critical exponent is $\gamma:=\gamma(\mu)=\inf_{[0,\infty)}\mbox{ supp}\, \tau_\circ=\inf_{S^{d-1}} \theta(u)$, e.g. if $d=1$, $\gamma=\inf \{\theta(1), \theta(-1)  \}$. Furthermore according to Corollary \ref{cor:testCE}, CTS$_\alpha$ laws are in $\mathcal S(\gamma)$ since their Thorin measures are absolutely continuous, and hence have no atoms. A cumulant analysis is again possible  using Proposition \ref{prop:momprop}, the integrations producing the same expression as \cite{kuc+tap:08}, Equation (2.11).  
\end{esmp}

\begin{esmp}\label{es:GRGTS}\emph{The Thorin measure of  generalized radially geometric stable laws}. Motivated by seeking tractable probability distributions with heavy tails but finite variance, in \cite{tor:24b} this author introduces the class of generalized radially geometric stable laws  GRGTS$_\alpha(\R^d)$, $\alpha \in (0,2)$. When $\alpha \neq 0$ we have  GRGTS$_\alpha(\R^d) \subset$ TS$_{\alpha}(\R^d)$ letting 
 \begin{equation}\label{eq:GRGTSlm}
 q(r,u)= e^{-\theta(u) r}E_{\rho(u)}\left(-\beta(u) r^{\rho(u)}\right)
 \end{equation}
in \eqref{eq:kts}, with $ \beta, \theta :S^{d-1} \rightarrow \R_+$, $\rho:S^{d-1} \rightarrow (0,1]$, $\rho, \beta$ bounded away from zero, and where $E_a$ is the Mittag-Leffler function defined as
\begin{equation}\label{eq:ML}
E_a(z)=\sum_{k=0}^\infty\frac{z^{k}}{\Gamma(\rho k+1)}, \qquad z \in \mathbb C, \quad a \in (0,1].
\end{equation} 
 It is easily shown that $q(\cdot, u)$ is c.m. for all $u \in S^{d-1}$. 
The obvious multivariate extension of \citet{tor:24b}, Proposition 3, to $d > 1$ states that 
\begin{align}\label{eq:sdensity}
Q_{u}( ds )&= \frac{(s-\theta(u))^{\rho(u)-1}}{\pi }\frac{ \sin( \rho(u) \pi )}{ \beta(u)^{-1} (s-\theta(u))^{2\rho(u)}+ 2 (s-\theta(u))^{\rho(u)} \cos(\rho(u) \pi)+\beta(u)}\I_{\{s \geq \theta(u)\}} ds.
\end{align}
  As a specific example, let us again consider an positively-supported case, i.e. $d=1$ and $Q_{-
 }=0$, and take $\alpha=1$ with $\rho, \theta, \beta$  some given constants. Let us try to compute the Thorin density $\tau$. In such case the fractional integral reduces to the usual Riemann integral and it is possible to perform the integration \eqref{eq:fracintmeas} explicitly.
Assuming for simplicity $\theta=0$,  $\beta=1$. The density in the right hand side reduces to
 \begin{equation}\label{eq:poldens}
s(x)=\frac{x^{\rho-1}}{\pi }\frac{ \sin(\rho \pi  )}{ x^{2\rho}+ 2 x^{\rho} \cos(\rho \pi)+1}\I_{\{x>0\}}
\end{equation}
 in which we recognize  the density of the Lamperti law (\citet{jam:10}),  which can be represented as the ratio $R=X_\rho/X'_\rho$ of two unit scale i.i.d. $X_\rho$, $X'_\rho$  r.v.s with law in S$_{\rho}(\R_+)$. The associated cumulative distribution function (c.d.f.) is 
\begin{equation}\label{eq:polcdf}
F_R(x):=1-\frac{1}{\rho \pi }\cot^{-1}\left(\cot(\rho \pi )+ \frac{x^{\rho}}{\sin(\rho \pi )}\right)\I_{\{x \geq 0\}},
\end{equation}
(\citet{jam:10},  Proposition 2.3).
Therefore we are led to the Thorin measure 
\begin{equation}
\tau=I_0^1 s(x) dx =F_R(x)dx.
\end{equation}
  The cases $\theta>0,\beta \neq 1$ are analogous and follow by a simple substitution in \eqref{eq:poldens} and \eqref{eq:polcdf} leading to the Thorin density $\tau(x)=F_R( (x-\theta)\beta^{-1/\rho} ) $.
 
\end{esmp}

It is also possible to introduce the  RGTS$(\R^d)$ class of radially geometric tempered stable distributions as the $T(\R^d)$ subclass having $k(u,r)=q(u,r)$, with $q$ in \eqref{eq:GRGTSlm}, whose associated Thorin family is given by $\{Q_u, \, u \in S^{d-1}\}$ with $Q_u$ as in \eqref{eq:sdensity}. The RGTS  class  is hence not a TS$_\alpha$ subclass. 
Such class is more commonly known as the tempered positive Linnik family; it was introduced in \citet{bar+al:16a}, and its associated processes were studied in \citet{tor+al:21}. Laws of this kind extend the geometric stable ones introduced in \citet{kle+al:84}, whose most famous instance is probably the \citet{pil:90} Mittag-Leffler distribution.


\subsection{Logistic type laws}\label{subsec:log}

An important class of distributions in applied social sciences and economics are the logistic laws, although perhaps their infinitely-divisible structure, until recently, has not drawn a lot of attention. For our purposes these offer an example of a non-gamma distribution class with  purely atomic Thorin measure.

Inspired by  \citet{gri:01}, a multivariate generalized logistic class GZD$(\R^d) \subset $ SD$(\R^d)$  can be   introduced  by letting
\begin{equation}\label{eq:GZDk}
 k(u,r)= \frac{e^{-c(u) r /\sigma(u)} }{1-e^{-r /\sigma(u)}} , \qquad  r>0
\end{equation}
with $c,\sigma : S^{d-1} \rightarrow \mathbb R_+$ continuous and bounded away from zero.
The case $d=1$, $c(-1)=c_-$, $c(1)=c_+$, $\lambda(du)=\lambda \delta_1(du)+\lambda \delta_{-1}(du)$, $\lambda>0$, $\sigma(\pm 1)=\sigma$ reduces to the univariate generalized logistic class introduced in  \citet{gri:01}, which we denote GZD$_G$. An important further subcase is the Meixner distribution, given by $c_\pm=1/2 \pm c$, $|c| < 1/2$ (\citet{gri:01}, \citet{sch:01}). The standard logistic and skew-logistic distributions are also clear special cases.

The  GZD$_G$ laws have an analytical c.f. in a neighborhood of 0  given by
\begin{equation}\label{gzCF}
\hat \mu(z)=\left(\frac{B\left(c_-+ i z \sigma , c_+-i z \sigma  \right)}{B(c_-, c_+)} \right)^{ \lambda}, \qquad z,  c_+, c_-, \lambda, \sigma >0.
\end{equation}
where $B(\cdot, \cdot)$ is Euler's Beta function.
When  $\lambda=1$, then $\mu$ belongs to  the $z$-distribution class ZD$(\R)$, more commonly known as Generalized logistic type IV.  See also \citet{bal:13} and \citet{bn+al:82}. It is shown in \citet{bon:12}, Example 7.2.4, that ZD$(\R) \subset T(\R)$. This inclusion can be  extended to  GZD$(\R^d) \subset T(\R^d)$ and the Thorin measures can be explicitly determined.

\begin{theoremEnd}[proof at the end,
	no link to proof]{prop}
We have  {\upshape GZD}$(\R^d) \subset T(\R^d)$ and $\mu \in$ {\upshape GZD}$(\R^d)$ has Thorin family $\{\tau_u , \, u \in S^{d-1} \}$ given by
\begin{equation}\label{eq:tauseruies}
\tau_{u}=\sum_{k=0}^\infty \delta_{(k+c(u))/\sigma(u)}.
\end{equation}
\end{theoremEnd}
	\begin{proofE}
Consider the series expansion  
\begin{equation}\label{eq:serieslog}
\frac{x^{a}}{1-x^{b}}=\sum_{k=0}^{\infty} x^{a +k b}, \qquad a,  b \neq 0, \quad |x|<1
\end{equation} 
and use  \eqref{eq:serieslog} in \eqref{eq:GZDk} with $x=e^{-y}$, $y>0$, $a=c(u) r /\sigma(u)$, $b= r /\sigma(u)$. Since $e^{-c y }=\mathscr L(\delta_c, y)$, $c \in \R$, we substitute the Laplace integral in each of the terms of the series thus obtained and finally interchanging summation and integration produces \eqref{eq:tauseruies}. 
\end{proofE}
The Thorin density of logistic laws is hence purely atomic, so by Proposition \ref{prop:existsmomThor}, Theorem \ref{thm:subexp} and Corollary \ref{cor:testCE}, convolution equivalence and existence of the absolute exponential moment at the infimum of the Thorin measure support do not hold.

\begin{esmp}\emph{The  } GZD\emph{$_G(\R)$ class}.
When $d=1$,  we show that performing the integration \eqref{eq:harexpThor} with  Thorin measures in \eqref{eq:tauseruies}, we recover \eqref{gzCF}. With $z \in \R$, we use the truncation function 1, since the moments exist, obtaining
\begin{align}
\hat \mu(z)=&\exp\left(- \left( \int_{0}^{\infty}  \left(\log \left(1  - \frac{iz}{s} \right)  + \frac{iz}{s}  \right) \tau_+(ds)+ \int_{0}^{\infty} \left( \log \left(1  + \frac{iz}{s}  \right)- \frac{iz}{s}\right)  \tau_{-}(ds) \right) \right)= \nonumber \\ =&\exp\left(-\lambda \sum_{k=0}^\infty \log \left[\left( 1 + \frac{iz  \sigma}{k+c_-} \right)  \left(1 - \frac{iz  \sigma}{k+c_+}  \right)  \right]  + iz \lambda \sigma  \left( \sum_{k=0}^\infty \frac{1}{(k+c_-)}- \sum_{k=0}^\infty \frac{1}{(k+c_+)} \right) \right) \nonumber \\=&\prod_{k=0}^{\infty}\left( 1 + \frac{iz \sigma}{k+c_-} \right)^{-\lambda}  \left(1 - \frac{iz \sigma}{k+c_+}  \right)^{-\lambda} e^{ iz \lambda  \sigma(\psi^{(0)}(c_+)-\psi^{(0)}(c_-))} \nonumber \\ =&\left(\frac{B\left(c_-+ i z \sigma , c_+-i z \sigma  \right)}{B(c_-, c_+)} \right)^{ \lambda}  e^{ -iz \lambda \sigma(\psi^{(0)}(c_-)-\psi^{(0)}(c_+))}
\end{align}
having used  \citet{gra+ryz:14}, 8.325.1 in the last equality and the series representation of the Digamma function $\psi^{(0)}(z):=\frac{d}{dz}\log \Gamma(z)$, $z \in \mathbb C \setminus \mathbb Z_{\leq 0}$ (\citet{abr+ste:48}, 6.3.16). 

To corroborate the line of analysis of subsections \ref{subsec:em} and \ref{subsec:ud}, notice that the m.g.f. of $\mu$ with 0 location is 
\begin{equation}\label{eq:mugzdg}
\hat \mu(-iz)=\left(\frac{B\left(c_- + z \sigma , c_+ - z \sigma  \right)}{B(c_-, c_+)} \right)^{ \lambda}.
\end{equation} An absolute exponential moment exists finite if and only if so do both the corresponding exponential moments of $\mu$ and $\check \mu$. But now  one has indeed $\gamma(\mu)=\inf_{S^{d-1}}\mbox{supp}\, \tau_\circ=\inf\{c_+/\sigma, c_-/\sigma\}$ from \eqref{eq:tauseruies}. If $c_+ \leq c_-$ then $\hat \mu(-i \gamma(\mu))=\hat \mu(-i c_+/\sigma)=\infty$ because the Beta integral only converges for positive value of the parameters. The same conclusion holds if $c_-<c_+$ since  $\check \mu$ is  obtained by interchanging $c_+$ and $c_-$ in \eqref{eq:mugzdg}. Therefore we verified the m.g.f. at $\gamma(\mu)$ with $\mu \in $ GZD$_G(\R)$ does not exists. Since the support of GZD distributions is purely discrete, this is in line with   Corollary \ref{cor:testCE}.
\end{esmp}


\begin{esmp}\emph{Gumbel distribution}.
The univariate Gumbel extreme value distribution $\mu$ with c.d.f. $F(x)=\exp\left( - e^{-\frac{x-c}{\sigma}} \right)$, $x,c \in \R, \sigma>0$ is related to the logistic family, since a standard logistic law can be obtained as a convolution of a Gumbel and a reflected Gumbel distribution. The class of Gumbel laws $G(\R)$ is such that $G(\R) \subset $ SD$(\R)$ although using standard methods this is not entirely straightforward to see.
We  show here that $G(\R) \subset$ GZD$(\R) \setminus $GZD$_G(\R)$.  
Consider the Thorin measure $\tau$ in \eqref{eq:tauseruies} to be positively-supported,  with constant function parameters $c(1)=1$, $\sigma(1)=\sigma>0$, but $\lambda(du)= \delta_1$. We obtain 
\begin{align}
\hat \mu(z)=& \exp \left(-  \int_{0}^{\infty}  \left(\log \left(1  - \frac{iz}{s} \right)  + \frac{iz}{s}  \right) \tau(ds)\right) \nonumber \\ =& \exp\left(-\sum_{k=0}^\infty \left(\log \left( 1 - \frac{iz  \sigma}{k+1}\right) +\frac{iz \sigma }{k+1}\right)  \right) \nonumber  \\ =& \exp \left(\log \left(\Gamma(-iz\sigma )(-iz \sigma) \right) - i z \sigma \gamma \right)  \nonumber \\ =& \Gamma(1-i z \sigma) e^{-i z \gamma \sigma} .
\end{align}
We have used the series representation of $\log \Gamma(z)$ in \citet{bor+mol:04} p.204, and here $\gamma$ stands for the Euler-Mascheroni constant. 
The last line above  is the c.f. of a centered $G(\R)$ distribution (\citet{ste+vh:03}, Example 9.15). 
As a result, elements of $G(\R)$ satisfy \eqref{eq:GZDk} but not \eqref{gzCF}. The Gumbel class is a  spectrally positive distribution class supported on $\R$, hence its associated L\'evy processes are necessarily of infinite variation. Again, convolution equivalence does not yield since $\gamma(\mu)=1/\sigma$ and the gamma function has a pole in 0. This is again in accordance with Corollary \ref{cor:testCE} and the discrete nature of the support of the Thorin measure. 

It is  thus possible to define a seemingly novel $d$-dimensional Gumbel class $G(\R^d)$ using \eqref{eq:GZDk} and letting $\lambda$ be supported on the positive orthant.
\end{esmp}

\subsection{A counterexample: the Dickman distribution}

A significant example of a distribution class contained in SD$(\R_+) \setminus T(\R_+)$  is the so-called Dickman distribution class $\mathcal D(\R_+)$, based on a c.d.f. belonging to a class of functions introduced in \citet{dic:30}. Its associated subordinators have been recently studied in  \citet{gup+al:24}; see references therein for further details. A distribution $\mu \in \mathcal D(\R_+)$ is defined through the canonical function
\begin{equation}
k(x)= \theta \I_{\{x  \in (0,1)\}}, \qquad x, \theta>0.
\end{equation}
Hence $k$ is not a Laplace transform, implying $\mathcal D(\R_+)\subset $ SD$(\R_+) \setminus T(\R_+)$. What one also observes is that being the L\'evy measure of the Dickman distribution compactly-supported, then $\int_1^\infty e^{\theta x}k(x) x^{-1} dx  <\infty$ for all $\theta \in \R$, and hence, being finiteness of truncated exponential L\'evy moments  equivalent to that of the distributional exponential moments, the bilateral Laplace transform of the Dickman law has abscissae of convergence $\pm \infty$. As a consequence, the c.f. of the Dickman distribution extends to an entire function, a property which is not shared by non-Gaussian Thorin distributions in light  of Corollary \ref{cor:analchar}.

\section{Subordination  theory for Thorin processes}\label{sec:sub}
Thorin processes are particularly well-behaved with respect to L\'evy subordination. Many  well-known formulae find an alternative simple interpretation when subordination is seen through the lenses of the Thorin representation.

Let us briefly recall that for a L\'evy process $X=(X_{
t})_{t \geq 0}$ on $\R^d$ and a subordinator $T$ i.e. an a.s. increasing L\'evy process $T=(T_t)_{t \geq 0}$,  the subordination of $X$ by $T$ is the stochastic process $Y=(X_{T_t})_{t \geq 0}$. It can be easily shown  that $Y$ is a L\'evy process, whose characteristic exponent is obtained as the functional composition of the Laplace exponent of $T$ with the characteristic exponent of $X$. The L\'evy triplet of $Y$ is also fully identifiable; see \citet{sat:99}, Theorem 30.1. In particular if $\rho$ is the L\'evy measure of $T$, and letting $\mu_t=\mathcal L(X_t)$ the L\'evy measure $\nu'$ of $Y$ is given by the Bochner integral
\begin{equation}\label{eq:bochint}
\nu'(B)=\int_{(0,\infty)}\mu_s(B)\rho(ds), \qquad  B \in \mathcal B(\R^d \setminus \{0 \}).
\end{equation}
For a fully detailed overview see \citet{sat:99}.

\subsection{Thorin subordination and potential measures}

We begin from a general result on the relationship between the L\'evy measure of the subordinated process and the Thorin measure of the subordinator. The two are connected through the L\'evy measures of the potentials of the parent L\'evy process.

 Recall that the $q$-th potential measure of a L\'evy process $X=(X_t)_{t \geq 0}$ is defined as

\begin{equation}\label{eq:potmeas}
V^q(B)=\int_0^\infty  e^{-q t } \mu_t(B) dt ,\qquad B \in \mathcal  B(\R^d \setminus \{0 \}).
\end{equation}
with $\mu_t=\mathcal L(X_t)$.
Notice that for all $q>0$, $q V^q$ is a probability measure and the induced distribution is i.d. with triplet $(0,0, \nu_q^*)$ and
\begin{equation}\label{eq:Levypotmeas}
\nu_q^*(B)=\int_0^\infty \frac{ e^{-q t }}{t} \mu_t(B) dt , \qquad B \in \mathcal  B(\R^d \setminus \{0 \}).
\end{equation}
See \citet{sat:99}, Definition 30.9 and Theorem 30.10.
When $T$ is a Thorin subordinator, the L\'evy measure of $Y=(X_{T_t})_{t \geq 0}$ can be expressed by means of a Bochner-type integral of the measures $\nu_q^*$ with respect to the Thorin measure of $T$.


\begin{theoremEnd}[proof at the end,
	no link to proof]{thm}\label{thm:ThorinSubBMSym}
Let $X=(X_t)_{t \geq 0}$ be a L\'evy process on $\R^d$ with L\'evy triplet $(a, \Sigma, \nu)$
and $T=(T_t)_{ t \geq 0}$ a Thorin subordinator with Thorin triplet $(b, 0, \tau)$. Then $Y=(X_{T_t})_{t\geq 0}$ is a L\'evy process with L\'evy  triplet $(a', \Sigma', \nu')$ with 
\begin{align}
a' &= a \, b + \int_{\{|x| \leq 1\}}\int_{0}^\infty x \nu^*_s(dx) \tau(ds) \label{eq:ThorinLevySubdrift} \\
\Sigma'&= b \Sigma   \label{eq:ThorinLevySubQuad}\\ 
\nu'(B)&= b \, \nu(B)+ \int_{0}^\infty \nu^*_s(B) \tau(ds), \qquad B \in \mathcal B(\R^d \setminus \{0 \}) \label{eq:ThorinLevySubLM}
\end{align}
where $\nu^*_s$ is the L\'evy measure of $sV^s$.
\end{theoremEnd}

\begin{proofE} In view of \citet{sat:99}, Theorem 30.1, \eqref{eq:ThorinLevySubQuad} follows directly. Looking at \eqref{eq:ThorinLevySubLM}, the genesis of the first summand is also clear from such theorem. Letting then $\rho$ be the L\'evy measure of $T$, denoting $\mu_t=\mathcal L(X_t)$ and using  \eqref{eq:Levypotmeas}, after interchanging integration order by Fubini's Theorem, we have that
\begin{align} 
\int_0^\infty \mu_t(B) \rho(dt)=\int_{0}^\infty \int_0^\infty \frac{e^{-t s}}{s} \mu_t(B) \tau(ds)  =\int_0^\infty \nu^*_s(B) \tau(ds) ,\quad B \in \mathcal B(\R^d)
\end{align}
and \eqref{eq:ThorinLevySubLM} follows again from \citet{sat:91}, Theorem 30.10, and so does Equation \eqref{eq:ThorinLevySubdrift}.
 \end{proofE}
A noticeable application is that to Brownian subordination with respect to gamma processes.

\begin{esmp}\emph{Symmetric variance gamma  process}. Assume $d=1$,  that  $T$ is a gamma subordinator with $T_1 \in \Gamma(\theta, \theta)$, Var$[T_1]=\theta^{-1}$, and $X$ is a standard Brownian motion. Then $Y=X_T$ is the symmetric version of the so-called  \emph{variance gamma process} (VG, \citet{mad+al:98}). The measure $\nu_q^*$ determined by $X$ is absolutely continuous, and its density can be obtained from \citet{sat:99}, Equation (30.30). With $K_p(x)$ the modified Bessel function of the second kind with order $p$ we have, recalling $K_{1/2}(x)=e^{-x}/\sqrt{x}\cdot \sqrt{\pi/2}$, that the following holds
\begin{equation}
\nu^*_q(x)= \sqrt{\frac{2}{\pi}}\frac{(2 q)^{1/4}}{\sqrt{|x|}}K_{1/2}(\sqrt{2 q}|x|)=   \frac{e^{-\sqrt{2 q}|x|}}{|x|}.
\end{equation}
Applying \eqref{eq:ThorinLevySubLM} with the Thorin measure $\tau=\theta \delta_\theta$ of $T$, yields the following well-known  (e.g. \citet{mad+al:98})
 expression  for the L\'evy density  $\nu'$ of $Y$ 
\begin{equation}
\nu'(x)=\theta \frac{e^{-\sqrt{2 \theta}|x|}}{|x|}.
\end{equation}
We will show how to determine the L\'evy measure from the Thorin measure of a general VG  process in the next section.
\end{esmp}

\subsection{Brownian subordination}\label{subsec:BrownSub}


We investigate next how univariate drifted  Brownian motions interact with GGC subordinators, and show that such subordination produces processes in the Thorin class. Moreover, we explicitly determine the Thorin triplet of the subordinated processes. Conversely, we establish necessary and sufficient conditions for a Thorin process to be interpreted as a GGC-subordinated drifted Brownian motion, and find the Thorin measure of such GGC subordinator.

 It is known from \citet{bon:12}, Theorem 7.3.1, that a standard Brownian motion on $\R$ subordinated to a GGC process is distributed as a Thorin process. We extend this result to Brownian motions with drift and explicitly determine the Thorin triplet. This generalization is non-trivial and important, as it makes possible to construct  distributionally asymmetric processes via subordination.

\begin{theoremEnd}[proof at the end,
	no link to proof]{thm}\label{thm:ThorinSubBM1} Let $X=(X_t)_{t \geq 0}$  be given by $X_t=\sigma W_t+\theta t$,  $\sigma >0$, $\theta \in\R$, with $W=(W_t)_{t \geq 0}$ a standard Brownian motion on $\R$, and let $T=(T_t)_{t \geq 0}$ be a  Thorin (GGC) subordinator independent of $X$ with Thorin triplet $(a, 0, \rho)$ . Then $Y=(X_{T_t})_{t \geq 0}$ is a Thorin process having Thorin triplet $(b, a \sigma^2, \tau)$ with
\begin{align}
b&=  \theta \left(a  + \int_{ \{t \geq \sigma^2/2+|\theta|\}}\frac{ \rho(dt)}{t}  -\frac{ \mbox{\upshape{sgn}}(\theta)}{\theta} \int_{ [\sigma^2/2-|\theta|,\sigma^2/2+|\theta| ) }\frac{ \rho(dt)}{s^{-|\theta|/\sigma^2}\circ f(t)} \right)    \label{eq:driftTC}\\ 
\tau_+ &=  (f_\sharp \rho)^{\theta/\sigma^2}, \qquad \tau_-=  (f_\sharp \rho)^{-\theta/\sigma^2},    \label{eq:tmTC}
\end{align}
 where $f:[0, \infty) \rightarrow  \left[ \frac{|\theta|}{2 \sigma^2}, \infty \right) $ is the hyperbolic function 
 \begin{equation}\label{eq:f}
 f(x)= \frac{\sqrt{\theta^2+ 2 \sigma^2 x}}{\sigma^2}.
 \end{equation}
 Furthermore,  if $a=0$ and for all $p\in [0,2]$, we have $\beta(T) \leq p/2$ if and only if $\beta(Y) \leq p$. Likewise, $E[|Y_t|^p]<\infty$ if and only  $E[T_t^{p-\frac{p}{2} \I_{\{\theta = 0\}}} ]<\infty$.  

\end{theoremEnd}
\begin{proofE}

We know from general theory (\citet{sat:99}, Theorem 30.1) that $Y$ is a L\'evy process. With  $\mu=\mathcal L(T_1)$,  by conditioning under independence and denoting $ n $ a normal $N(\theta, \sigma^2)$ law, we have the familiar relationship
\begin{align}
\psi_Y(z)&=\log E[e^{i z Y_1}] =\log E[ e^ { i  T_1 \log \hat n(z)}]=\log E[ e^{  i T_1 (-\sigma^2z^2/2+ i \theta  z ) }]  \nonumber \\ &= \log \hat \mu(-\sigma^2 z^2/2+ i \theta  z) =   i a \theta z- a \sigma^2  z^2/2  + 
\int_{0}^\infty \log\left( \frac{t}{t - i z \theta +  \sigma^2 z^2/2} \right) \rho(dt) .
\end{align}
 Let $f^\pm:[0,\infty)\rightarrow \R_+$, $f^\pm(x)= \mp \theta/\sigma^2 + \sqrt{\theta^2 + 2 \sigma^2 x}/\sigma^2$, $x \geq 0$, that is  $f^+=s^{\theta/\sigma^2} \circ f$, $f^-=s^{-\theta/\sigma^2} \circ f$.  
First we show that $\tau$ corresponding to the family $\{\tau_+, \tau_-\}$  satisfies \eqref{eq:thorconv0}.  Observe $\{ f^+(x) \geq 1 \}=\{ x \geq \theta+\sigma^2/2\}$, $\{ f^-(x) \geq 1 \}=\{ x \geq \theta-\sigma^2/2\}$, and recall that $\int_1^\infty t^{-1} \rho(dt)<\infty$ since a subordinator must necessarily be a process of finite variation. Therefore
\begin{align}
\int_{\{y \geq 1\}} \frac{ (f_\sharp\rho)^{\pm \theta/\sigma^2}(dy)}{y^2}
=\int_{\{ t \geq \sigma^2/2 \pm \theta \}} &\frac{\sigma^4 \rho(dt)}{2 \theta^2 \mp 2 \theta \sqrt{\theta^2+2 \sigma^2t } + 2 \sigma^2 t} <  \infty
\end{align}
 follows for all $\theta \in \R$. 
Furthermore it is
\begin{align}
\int_{\{y \leq 1\}}  |\log  y |(f_\sharp\rho)^{\pm \theta/\sigma^2 }(dy)=\int_{\{0 \leq t\leq   \sigma^2/2 \pm \theta \}}&  \left|\log \left(\frac{\mp \theta   + \sqrt{\theta^2+ 2 \sigma^2 t}}{\sigma^2} \right)\right| \rho(dt) <\infty 
\end{align}
as a consequence of \eqref{eq:thorconv0} holding for $\rho$. This shows that the integrability conditions  for $\tau$ in  \eqref{eq:tau1d} hold when $\tau_+,\tau_-$ are as in  \eqref{eq:tmTC}. Next we need to show the Thorin--Bondesson representation. It holds
\begin{align}\label{eq:nodrifted}
&\int_{0}^\infty  \log\left( \frac{t}{t - i z \theta +  \sigma^2 z^2/2} \right) \rho(dt) + i a \theta z&  \nonumber \\ &= \int_{0}^\infty  \left( \log\left( \frac{-\theta + \sqrt{\theta^2 + 2 \sigma^2 t}}{-\theta + \sqrt{\theta^2 + 2 \sigma^2 t} - i \sigma^2 z} \right) +  \log\left( \frac{-\theta - \sqrt{\theta^2 + 2 \sigma^2 t}}{-\theta - \sqrt{\theta^2 + 2 \sigma^2 t} - i \sigma^2 z} \right)  \right)  \rho(dt)   + i a \theta z  \nonumber \\ &= \int_{0}^\infty  \left( \log \left( \frac{f^+(t)}{f^+(t) - i z} \right)  +   \log \left( \frac{f^-(t)}{f^-(t) + i z} \right)   \right) \rho(dt) + i a \theta z.
\end{align}
Now we have
\begin{equation}\label{eq:f+}
\{ f^{\mbox{sgn}(\theta)}(x) \geq 1\}=\left\{ x \geq  \frac{\sigma^2}{2}+|\theta|\right\},
\end{equation}
and
\begin{align}\label{eq:f-}\{ f^{\mbox{-sgn}(\theta)}(x)\geq 1\}&=\left\{ x \geq  \frac{\sigma^2}{2}-|\theta|\right\}= \left\{x \geq \frac{\sigma^2}{2}+|\theta|\right\} \bigcup\left[\frac{\sigma^2}{2}-|\theta|, \frac{\sigma^2}{2}+|\theta|\right) \nonumber \\ &=\{ f^{\mbox{sgn}(\theta)}(x) \geq 1\} \bigcup\left[\frac{\sigma^2}{2}-|\theta|, \frac{\sigma^2}{2}+|\theta|\right)
\end{align} 
where the union is disjoint. From \eqref{eq:f-},  \eqref{eq:f+}, the observation that $1/f^+(x)-1/f^-(x)=\theta x^{-1}$, $x>0$, together with $f^{-\mbox{sgn}(\theta)}= s^{-|\theta|/\sigma^2} \circ f$,   we can  rearrange the integral \eqref{eq:nodrifted} as follows
\begin{align}\label{eq:comp}
&\int_{0}^\infty  \left( \log \left( \frac{f^+(t)}{f^+(t) - i z} \right)  +   \log \left( \frac{f^-(t)}{f^-(t) + i z} \right)   \right) \rho(dt) + i a \theta z
  \nonumber 
\\ &= \int_{0}^\infty  \left( \log \left( \frac{f^+(t)}{f^+(t) - i z} \right)  +   \log \left( \frac{f^-(t)}{f^-(t) + i z} \right) - \left(\frac{iz}{f^+(t) }\I_{\{f^+(t)\geq 1\}}-\frac{iz}{f^-(t) }\I_{\{f^{-}(t)\geq 1\}} \right)  \right) \rho(dt) + \nonumber \\ &\phantom{xxxxxxxxxx}+ia \theta z + i z\left(\int_{ \{ f^{\mbox{{\footnotesize sgn}}(\theta)}(x) \geq 1\}} \frac{\theta}{t}\rho(dt)  - \mbox{\upshape{sgn}}(\theta) \int_{ [\sigma^2/2-|\theta|,\sigma^2/2+|\theta| ) }  \frac{\rho(dt)}{f^{-\mbox{sgn}(\theta)}(t) } \right) \nonumber \\
&=\int_{-\infty}^\infty  \left( \log \left( \frac{y}{y - i z} \right)  -  \frac{iz}{y }\I_{\{|y| \geq 1 \}   } \right) \tau(dt) + i b \theta z  
\end{align}
which simultaneously shows \eqref{eq:driftTC} and \eqref{eq:tmTC}. It is critical to notice the convergence of the second integral term in the third line \eqref{eq:comp}, again as an effect of $\rho$ being the Thorin measure of a GGC process.

For the second to last statement we require  $a=0$, else a Brownian motion is present in the process. We observe, along the same lines as above, that \begin{align}
\int_{\{|y| \geq 1\}} \frac{\tau(dy)}{|y|^p}=\int_{\{t \geq \sigma^2/2  +\theta\}} &\frac{\sigma^2 \rho(dt)}{(-\theta+\sqrt{2 \sigma^2 t+ \theta^2})^p } +\int_{\{t \geq \sigma^2/2  -\theta \}} \frac{\sigma^2 \rho(dt)}{(\theta+\sqrt{2 \sigma^2 t+ \theta^2})^p}
\end{align}
where the finiteness of either side follows if and only if $\int_1^\infty t^{-p/2}\rho(dt)<\infty$ which by Proposition \ref{prop:fv}, is equivalent to  $\beta(T)\leq p/2$. Regarding expectations, the test integral for finiteness of the $p$-th absolute moment of $Y_t$, by Proposition \ref{prop:momprop}, because of the relationships complementary to \eqref{eq:f+} and \eqref{eq:f-}, writes as
\begin{align}
\int_{\{|y| < 1\}} \frac{\tau(dy)}{|y|^p}&=\int_{( 0,\sigma^2/2  -|\theta|)} \left(\frac{1}{(f^+(t))^p}+\frac{1}{(f^-(t))^p} \right)\rho(dt)+ C \nonumber \\ &=\int_{(0, \sigma^2/2  -|\theta|)}  \frac{ \left(\theta +\sqrt{\theta^2+2 \sigma^2 t}\right)^p + \left(-\theta+\sqrt{\theta^2+2 \sigma^2 t}\right)^p}{  2^pt^p}  \rho(dt) +C
\end{align}
for some $C>0$.
If $\theta \neq 0$ (respectively, $\theta=0$), the integrals at both sides converge if and only if $\int_{(0,1)} \rho(dt)t^{-p}<\infty$ (respectively, $\int_{(0,1)} \rho(dt)t^{-p/2}<\infty$ ), which is equivalent to $E[T_t^{p}]<\infty$ (resp. $E[T_t^{p/2}]<\infty$), again by Proposition 3. This concludes the proof of the Theorem.
\end{proofE}

As a corollary we obtain  that the canonical density and L\'evy measure
 of the subordinated process,  
can be written, respectively, as an integral transform of the Thorin measure of the subordinator  and an expectation of its canonical function.

Recall that the inverse Gaussian IG$(a,b)$ law is the absolutely continuous probability distribution with p.d.f. $h$ given by
\begin{equation}\label{eq:IG}
h(x)=\sqrt{\frac{b}{2 \pi x^3}} e^{-\frac{b(x-a)^2}{2 a^2 x}}\I_{\{x>0\}}, \qquad a,b>0.
\end{equation}
The result is as follows.

\begin{theoremEnd}[proof at the end,
	no link to proof]{cor}\label{cor:cansub}
Under the assumptions of Theorem \ref{thm:ThorinSubBM1}, the canonical function $k$ of $Y$ is given by 
 \begin{equation}\label{eq:cansub}
 k(x)=e^{x \theta/\sigma^2} \int_0^\infty e^{-|x|\frac{ \sqrt{\theta^2 + 2 \sigma^2 s}  }{\sigma^2}} \rho(ds), \qquad x \neq 0 ,
 \end{equation}
 and moreover the L\'evy density $\nu(x)$ of $Y$ has the probabilistic representation
\begin{equation}\label{eq:levsub}
\nu(x)=\left\{ \begin{array}{cc }E[k_\rho( Z_x)] & \mbox{ if  {\upshape sgn}}(\theta)=\mbox{{\upshape sgn}}(x)  \\ e^{2 x \theta/ \sigma^2}E[k_\rho(Z_x)] & \mbox{ if  {\upshape sgn}}(\theta)\neq\mbox{{\upshape sgn}}(x) \end{array}\right.
\end{equation} 
 where $k_\rho$ is the canonical function of $T$ and  $\mathcal L(Z_x) \in$ {\upshape IG}$(|x/\theta| ,(x/\sigma)^2)$.
\end{theoremEnd}

\begin{proofE}Let $\tau$ be the Thorin measure of $Y$.
Using \eqref{eq:tmTC}, and \eqref{eq:thorinm}, for $x \neq 0$ and the notation in Theorem \ref{thm:ThorinSubBM1} we have
\begin{align}\label{eq:canonical}
k(x)=&\int_{0}^{\infty} e^{- s|x|}\tau(ds)= 
\left(\int_{0}^{\infty} e^{- x  f^+(s)}\rho(ds)\right) \I_{\{x > 0\}} + \left(\int_{0}^{\infty} e^{ x  f^-(s)}\rho(ds) \right) \I_{\{x < 0\}} 
\end{align}
and  this proves \eqref{eq:cansub}.
For $\eta  \in$ IG$(a,b)$ it is well-known that
\begin{align}\label{eq:IGlt}
\mathscr L(\eta, s)=\exp\left(\frac{b}{a}- \sqrt{b}\sqrt{b/a^2 + 2s } \right).
\end{align}
Choosing $a=x/\theta, b=(x/\sigma)^2$, whenever $x/\theta>0$ i.e $\mbox{sgn}(x)=\mbox{sgn}(\theta)$,
 \eqref{eq:IGlt} specifies to $\mathscr L(\eta, s)$, with $\eta \in$ IG$(x/\theta, (x/\sigma)^2)$. 	
If instead $\mbox{sgn}(x)\neq \mbox{sgn}(\theta)$, it coincides with $e^{2 x \theta/ \sigma^2} \mathscr L(\eta', s) $ with $\eta' \in$ IG$(-x/\theta, (x/\sigma)^2)$.
In either case, 
using Fubini's Theorem to interchange the order of integration, denoting with $p_x$ the density of $ Z_x$ and letting $q_x=p_x$ or  $q_x=e^{2 x \theta/ \sigma^2}p_x$ as required by the sign relationship of $\theta$ and $x$, we obtain
\begin{equation}\label{eq:cansub2}
 k(x)=e^{x \theta/\sigma^2} \int_0^\infty e^{-|x|\frac{ \sqrt{\theta^2 + 2 \sigma^2 s}  }{\sigma^2}} \rho(ds)=\int_0^\infty \rho(ds)\int_0^\infty e^{- u s } q_x(u)du
=\int_0^\infty k_\rho(u) q_x(u) du.
   \end{equation}
   which finishes the proof.
\end{proofE}

 Corollary \ref{cor:cansub} provides an alternative way to  \eqref{eq:bochint} of calculating L\'evy measures of Brownian motions subordinated to GGC subordinators. In the examples below we let $X$ be the drifted Brownian motion of Theorem \ref{thm:ThorinSubBM1} and  $Y=(X_{T_t})_{t \geq 0}$.

\begin{esmp}\label{eq:varg}\emph{Variance gamma.} Consider for 	$T$ a gamma $\Gamma(\lambda, \beta)$, leading to a VG process for $Y$. From Subsection \ref{subsec:gamma} and  Theorem \ref{thm:ThorinSubBM1} we immediately have that the Thorin family of $Y$ is given by
\begin{equation}
\tau_\pm(ds)=\lambda \delta_{\frac{-\theta \pm \sqrt{2 \sigma^2 \beta+\theta^2}}{\sigma^2}}(ds).
\end{equation}
 From \eqref{eq:cansub} or \eqref{eq:levsub} it is  straightforward to derive the well-known canonical function $k$ of the variance gamma process i.e.
\begin{equation}
k(x)=\lambda  \exp\left(x \frac{\theta}{\sigma^2} -|x| \frac{\sqrt{2 \sigma^2 \beta+\theta^2}}{\sigma^2}\right).
\end{equation}


Thorin measures and canonical functions of Brownian motions subordinated with respect to finite convolutions of gamma processes are obtained similarly. This example explains why  the VG canonical function has the same functional expression as the a Laplace transform of an IG law.
\end{esmp}

\begin{esmp}\label{es:nts}\emph{Normal tempered stable.}  A fairly general  Brownian subordinated class of processes can be generated by using as a subordinator $T$ a CTS$_\alpha(\R_+)$ process, i.e. a Thorin process $T$ with canonical function $k_\rho(x)=\lambda x^{-\alpha} e^{-\beta x}$, $\alpha \in (0,1)$, $x,\lambda, \beta >0$. The process $Y$  is called a normal tempered stable process (NTS). As shown in Subsection \ref{subsec:TS} Thorin measures in the CTS$_\alpha(\R^d)$ class are absolutely continuous. Using the change of variable entailed by \eqref{eq:tmTC} and the density in \eqref{eq:CTStm} with $d=1$ we obtain that the NTS Thorin density  is as follows
\begin{align}\label{eq:taunts}
\tau(y)=\frac{\lambda}{\Gamma(\alpha)} \Bigg((\sigma^2 y^2/2+\theta y -\beta)^{\alpha-1}( \sigma^2 y +\theta) &\I_{\left\{y \geq \frac{-\theta+\sqrt{\theta^2+ 2 \beta \sigma^2}}{\sigma^2}\right\}}  \nonumber \\ &  +(\sigma^2 y^2/2 -\theta y -\beta)^{\alpha-1}( \sigma^2 y -\theta)\I_{\left\{y \leq \frac{-\theta-\sqrt{\theta^2+ 2 \beta \sigma^2}}{\sigma^2}\right\}}  \Bigg). 
\end{align}
 The probabilistic representation \eqref{eq:levsub} is well suited to determine the NTS canonical function $k$.
 Assume $\theta \neq 0$. Let $\xi=(x/\sigma)^2$, $\epsilon=x/\theta$; then
we have
\begin{align}\label{eq:IGexp}
k(x)&=\lambda\Big(\I_{\mbox{sgn}(\epsilon)=1}+e^{2 x \theta/\sigma^2}\I_{\mbox{sgn}(\epsilon)=-1}\Big)E[e^{-\beta Z_x}(Z_x)^{-\alpha}] \nonumber \\ &=\lambda \sqrt{\frac{\xi}{2 \pi}} e^{x\theta/\sigma^2 }\int_0^\infty u^{-\alpha-3/2 } \exp\left(-\frac{(\xi+ 2 \beta \epsilon^2) u}{2 \epsilon^2} -\frac{ \xi}{2 u} \right) du.
\end{align}
The generalized inverse Gaussian p.d.f. $g$, defined for $a,b>0$, $c \in \R$ is given by
\begin{equation}\label{eq:GIGpdf}
g(x)= \frac{(a/b)^{c/2}}{2 K_{c}(\sqrt{a b})}x^{c-1}\exp \left( -\frac{a x}{2}- \frac{b}{2 x}\right) \I_{\{x>0\}}
\end{equation}
with $K_{c}(\cdot)$ the modified Bessel function of the second kind. Now  take $c=-\alpha-1/2$, $a=(\xi+ 2 \beta \epsilon^2)/\epsilon^2$, $b=\xi$. Using such parameters  we can exploit \eqref{eq:GIGpdf} to work out the integral in \eqref{eq:IGexp}. 
Notice $\sqrt{ab}=|x|\sqrt{\theta^2 + 2\beta \sigma^2}/\sigma^2$ and $\sqrt{a/b}=\sqrt{\theta^2 + 2\beta \sigma^2}/|x|$ and recall $K_{-c}(\cdot)=K_c(\cdot)$. Unwinding the substitutions we arrive to the known expression
 \begin{equation}\label{eq:knts}
k(x)=\frac{2 \lambda }{\sqrt{2 \pi } \sigma} (\theta^2 + 2\beta \sigma^2)^{\alpha/2+1/4}  \frac{e^{x \theta/\sigma^2  }}{|x|^{\alpha-1/2}}    K_{\alpha+1/2}\left( |x|  \frac{\sqrt{\theta^2+ 2 \beta \sigma^2 }}{ \sigma^2} \right)           
 \end{equation}
see e.g. \citet{con+tan:03}, Equation (4.24). A similar calculation, but simpler, shows \eqref{eq:knts} when $\theta=0$ in \eqref{eq:taunts}. An important special case is $\alpha=1/2$, when $T$ specifies to an inverse Gaussian process (\citet{bn:97}), with canonical function
 \begin{equation}
k(x)=\frac{2 \lambda }{\sqrt{2 \pi } } \frac{\kappa}{\sigma}  e^{x \theta/\sigma^2  }    K_{1}\left( |x|  \frac{\kappa}{ \sigma^2} \right)                
 \end{equation}
 where $\kappa= (\theta^2 + 2\beta \sigma^2)^{1/2}$.
 
 Observe how the Thorin measure of a Normal tempered stable process is a simple rational function, whereas its L\'evy measure involves a special function.
\end{esmp}


We provide next an example in which the Thorin measure of the subordinated process can be easily identified, whereas it is not immediately clear how to proceed for the L\'evy one.

\begin{esmp}\label{es:NL}\emph{Normal Linnik.} Consider 
a subordinator $T$ having as a unit time marginal a RGTS$(\R_+)$ law as in Subsection \ref{subsec:TS}, i.e. a tempered  Linnik subordinator. The resulting subordinated process $Y$ could be called a  Normal Linnik process (NL). 
The c.f. of one such process is easily found applying the usual conditioning argument to $\phi_T$ as in \citet{bar+al:16a}. 

In the context of Example \ref{es:GRGTS} with $d=1$, take $\rho(1)=\alpha \in (0,1)$, $\beta(1)=\kappa >0, \theta(1)=\beta >0 $, together with a constant spherical function at $\lambda>0$.
As in the previous example, combining  \eqref{eq:tmTC} and  \eqref{eq:sdensity} we can see that the Thorin densities  are given by
\begin{align}\label{eq:tauNL}
\tau_\pm(y)&=\lambda \frac{(\sigma^2 y^2/2\pm \theta y-\beta)^{\alpha-1}}{\pi  }\times \nonumber \\  & \phantom{xxxxxx} \frac{ \sin(\alpha \pi  )(\sigma^2 y \pm \theta)}{ (\sigma^2 y^2/2\pm \theta y-\beta)^{2\alpha}/\kappa+ 2 (\sigma^2 y^2\pm \theta y-\beta)^{\alpha} \cos(\alpha \pi)+\kappa}\I_{\left\{y \geq   \frac{\mp \theta+ \sqrt{\theta^2+ 2 \beta \sigma^2 }}{\sigma^2}\right\}}.
\end{align}
The computation of the L\'evy measure using \eqref{eq:bochint}, would  involve a series of special functions, in view of the Mittag-Leffler function \eqref{eq:ML} appearing in the canonical function of the subordinator. Moreover, interchangeability of expectation and the series is far from clear, because of the lack of uniformity of convergence on account of the singularities in $x=0$ of all the series summands.  
\end{esmp}
\bigskip

It is well-known (e.g. \citet{che+shi:02},  Theorem  3.17), that not every L\'evy process can be represented as a subordinated Brownian motion. However,  within the Thorin class, whenever one such representation is possible, Theorem  \ref{thm:ThorinSubBM1} has a converse.

\begin{theoremEnd}[proof at the end,
	no link to proof]{thm}\label{thm:ThorinSubBM2}
Let $Y$ be a Thorin process with triplet $(0,0,\tau)$. 
Then there exists a Brownian motion $W=(W_t)_{t \geq 0}$, a GGC subordinator $T=(T_t)_{t \geq 0}$ independent of $W$ with Thorin triplet $(0,0,\rho)$, and a constant $\theta \in \R$ such that $Y=^d W_T +\theta T$, if and only if $\tau_+^{-\theta}=\tau^{\theta}_-$ and $\inf_{[0,\infty]} \mbox{{\upshape supp}}\, \tau_\circ \geq |\theta|$.  In such case the Thorin measure of $T$ is given by
\begin{equation}\label{eq:invpf}
\rho=g_\sharp \tau^{-\theta}_+=g_\sharp \tau^{\theta}_-
\end{equation}
with
 $g:[|\theta|,\infty) \rightarrow \R_+$ defined as
 \begin{equation}
  g(x)=\frac{x^2-\theta^2}{2 }.
\end{equation}  
\end{theoremEnd}
\begin{proofE}
For the if part, when $\sigma=1$ the function $f$ in \eqref{eq:f} is injective with positive domain and range $[|\theta|,\infty)$. Since $f^{-1}=g$ , we see that \eqref{eq:invpf} yields a well defined inverse pushforward, being supp$\, \tau_\circ \subseteq [|\theta|, \infty)$ by assumption. Inverting the pushforward maps in \eqref{eq:tmTC} yields \eqref{eq:invpf} and that supp$\, \rho \subseteq [0,\infty)$ is now clear. 
 We need to show that $\rho$ is a Thorin measure of a GGC law by verifying \eqref{eq:GGCtm}. But this follows since
\begin{equation}
\int_1^\infty \frac{\rho(dt)}{t}=\int_{\{ y \geq \sqrt{\theta^2 +2}\}} \frac{\tau_+^{-\theta}(dy)}{g(y)} =\int_{\{y \geq - \theta +\sqrt{\theta^2 +2}\}} \frac{\tau_+(dy)}{y^2/2 + \theta y  } <\infty, 
\end{equation}
 because $\tau$ is a Thorin measure. By a similar integration
\begin{equation}
\int_0^1 |\log t| \rho(dt)=\int_{\{y \leq - \theta +\sqrt{\theta^2 +2}\}} |\log (y^2/2+\theta y)|\tau_+(dy)<\infty 
\end{equation}
since $|\log (y^2/2+\theta y)| \sim b|\log y +c| $, $y \rightarrow 0^+$ for some constants $b,c>0$.

In order to prove the converse assertion, let us observe first that $\inf_{[0,\infty]} $ supp$\, \tau_\circ  > |\theta|   $ is necessary for  $\rho$ to be a Thorin measure of a GGC law in the first place. Assume thus  $\tau_+^{-\theta}\neq  \tau^{\theta}_-$, which after Laplace transforming is equivalent to  $e^{x \theta}k(-x)\neq k(x)e^{-x \theta }$. Then from \citet{con+tan:03}, Theorem  4.3,  \emph{2}, $Y$ a self-decomposable process with one such canonical function $k$ is not representable as a L\'evy-subordinated Brownian motion, and hence in a particular as a GGC-subordinated one. 
\end{proofE}

A symmetric Thorin  process $X$ requires $\tau_+=\tau_-$; on the other hand, a symmetric Brownian subordinated process must have $\theta=0$. Therefore from Theorem \ref{thm:ThorinSubBM2} it follows that all symmetric Thorin processes admit a Brownian subordinated representation. 
This special case was already contained in \citet{bon:12}, Theorem 7.13. Also taking into account Proposition \ref{prop:existsmomThor}, we have the following corollary    to Theorem \ref{thm:ThorinSubBM2}.

\begin{cor} Let $X$ be a Thorin process. If $X$ is symmetric then it can be represented as a subordinated Brownian motion. If $X$ is not symmetric then a necessary condition for it to be representable as a subordinated Brownian motion plus a linear drift is that $\mathcal L(X_1)$ admits a m.g.f. 
\end{cor}

From Theorem \ref{thm:ThorinSubBM2}, unwinding the pushforward map $g \circ s^{-\theta}(y)=y^2/2 - \theta y$ and using the inverse Laplace representation of $k$, also descends the following  counterpart of Corollary \ref{cor:cansub}.

\begin{cor}\label{cor:consubconv}
In the context of Theorem \ref{thm:ThorinSubBM2} we have that the canonical function $k_\rho$ of the subordinator $T$ can be expressed as, for $x > 0$ 
\begin{equation}\label{eq:cangauss}
k_\rho(x)=e^{ x \theta^2/2} \int_0^\infty e^{-x s^2 /2} \tau_+^{-\theta}(ds)=\int_\theta^\infty e^{-x (s^2/2 - s \theta)} \tau_+(ds),
\end{equation}
and the analogous equalities hold for $\tau_-$.
\end{cor}

The assumption of Theorem \ref{thm:ThorinSubBM2} that the Thorin measure is a shift of a symmetric one is not as constraining as it may appear. It is indeed  met for most distributions used in applications, among which those of Section \ref{sec:ThorList}. The exponentials appearing in the canonical functions are key to this remark, as can be seen in the examples below.

\begin{esmp}\emph{Bilateral gamma processes as a VG processes.} Assume we have a bilateral gamma process $Y$ with $\mathcal L(Y_1)=$Bil$\Gamma(\lambda, \theta_+,\lambda, \theta_-)$. The symmetry assumption in Theorem \ref{thm:ThorinSubBM2} making Brownian representability viable is satisfied as it can be seen by using \begin{equation}
A=(\theta_--\theta_+)/2, \qquad B= (\theta_-+\theta_+)/2
\end{equation} and then observing that the canonical function of $Y$ writes as $k(x)=\lambda e^{A x- B|x|}$, implying $\tau_+^{-A}=\tau_-^{A}= \lambda \delta_B$. The Dirac delta integration to which \eqref{eq:cangauss} reduces, produces \begin{equation}k_\rho(x)=\lambda e^{(A^2-B^2)x/2}=\lambda e^{-\theta_+\theta_- x/2}.\end{equation} Therefore, the bilateral gamma process can be represented in law as the VG process with $Y_t= W_{T_t} + (\theta_--\theta_+) t/2$ where $\mathcal L(T_1) \in \Gamma(\lambda, \theta_+\theta_-/2)$. This  has been known and used at least since \citet{mad+al:98}. 
\end{esmp}

\begin{esmp}\emph{Subordinated representation of CTS}  \emph{processes}. \citet{mad+yor:08} through a series of propositions covering a substantial part of the article, explain how to represent a CTS$_\alpha(\R)$ process as a subordinated Brownian motion with drift. 
Based on the discussion so far, their main Proposition 2 
trivializes. First of all notice that the assumptions in Theorem \ref{thm:ThorinSubBM2} are met also in the case of CTS$_\alpha(\R)$ processes with symmetric spherical part. Namely, assuming $\lambda_+=\lambda_-$ and letting again $A=(\theta_--\theta_+)/2$, $B=(\theta_-+\theta_+)/2$, leads to $k(x)=x^{-\alpha} e^{A x- B|x|}$, so that using  the Thorin measure found in Subsection \ref{subsec:TS} we deduce the required density relationship  $\tau_+^{-A}(y)=\tau_-^{A}(y)= 
\lambda(y-B)^{\alpha-1}(\Gamma(\alpha))^{-1}\I_{\{y \geq B\}} $.  It then follows from Corollary \ref{cor:consubconv}, and with the substitution $w=(s-B) \sqrt{x}$, that
\begin{align}\label{eq:rhoCTS}
k_\rho(x)=&\lambda \frac{e^{ x A^2/2}}{\Gamma(\alpha)}\int_B^\infty e^{-\frac{s^{2}x}{2}} (s-B)^{\alpha-1} ds=\lambda \frac{e^{ (A^2-B^2) x/2}}{\Gamma(\alpha)x^{\alpha/2}} \int_0^\infty e^{-\frac{w^{2}}{2} - w B \sqrt{x}} w^{\alpha-1} dw \nonumber \\ =& \lambda \frac{e^{-\theta_+ \theta_-x/2}}{x^{\alpha/2}} H_{-\alpha}\left( \frac{\theta_+ + \theta_-}{2} \sqrt{x}\right)
\end{align}
where $H_a$, $a<0$, is the Hermite function, given  by
\begin{equation}
H_a(z)= \frac{1}{\Gamma(-a)}\int_0^\infty e^{-x^2/2-xz}x^{-a-1}dx, \qquad z>0.
\end{equation}
The corresponding L\'evy density coincides with that determined in \citet{mad+yor:08}, Proposition 2.

\end{esmp}

\begin{esmp}\emph{Subordinated representation of generalized-$z$ processes}. Let $\mu \in$ GZD$_G(\R)$ with L\'evy measure of the form \eqref{eq:GZDk}, i.e. $c(\pm 1)=c_\pm >0$, $\sigma(\pm)=\sigma>0$ and $\lambda(du)=\lambda(\delta_1(du)+ \delta_{-1}(du)) $, $\lambda>0$. We can rewrite again $k(x)=e^{A x/\sigma -B|x|/\sigma}/(1-e^{-|x|/\sigma})$, $A=(c_--c_+)/2$, $B=(c_-+c_+)/2$ so we are under the assumption of Theorem \ref{thm:ThorinSubBM2}. Based on \eqref{eq:tauseruies}  and \eqref{eq:cangauss}, integrating term by term the series results in
\begin{align}\label{eq:newZrep}
k_\rho(x)&=\lambda e^{ \frac{A^2 x}{2\sigma^2} }\sum_{k =0}^\infty \exp\left(-x  \frac{(k+B)^2}{2 \sigma^2}  \right)=\lambda\sum_{k =0}^\infty \exp\left(-x  \frac{k^2+ 2 B k + B^2-A^2}{ 2 \sigma^2}\right) \nonumber \\ &=\lambda \sum_{k=0}^\infty \exp\left(-x  \frac{(k+c_+)(k+c_-)}{ 2 \sigma^2}\right).
\end{align}
The subordinator is then of the form of an infinite GGC convolution of gamma processes $G^k=(G^k_t)_{t \geq 0}$, $k \geq 0$, $\mathcal L(G^k_1)\in \Gamma\left(\lambda,\frac{(k+c_+)(k+c_-)}{ 2 \sigma^2}\right)$ whose limit in law is not otherwise known. The symmetric case $c_+=c_-$, $\sigma=1/\pi$ was treated in \citet{bn+al:82} and \citet{pit+yor:03} in connection with some infinitely divisible laws arising from hyperbolic functions. The Meixner case was obtained by some slightly involved arguments in \citet{mad+yor:08}. 
Some partial results on the asymmetric case are also offered in \citet{car+tor:21}, but to the best of our knowledge \eqref{eq:newZrep} in its full generality is new.
\end{esmp}

\subsection{Subordination of gamma processes by generalized negative binomial convolutions}

In the following we illustrate some aspects of Thorin subordination based on bilateral gamma subordinands. The natural subordinators to univariate drifted Brownian motions are infinite activity  L\'evy processes, and we dedicated the previous section to show that if the subordinator is Thorin, the subordinated process must be also Thorin.  When the driving process is  bilateral gamma, i.e. the prototypical  BG$_0(\R) \cap T(\R)$ process, the natural subordinator is instead of finite activity with a specific marginal distribution, namely, negative binomial. Such subordinator is therefore not even self-decomposable. Nevertheless, a closure property similar to that of Brownian subordination  holds, in that subjecting a univariate bilateral gamma processes  to negative binomial subordination yields L\'evy processes in BG$_0(\R) \cap T(\R)$. This closure property is maintained when enlarging the family of the subordinators considered to the smallest one containing finite convolutions of negative binomial processes and their weak limits.

\medskip

Recall that the probability generating function (p.g.f.)  of a r.v. $X$ is defined as $\varphi_X(z):=E[z^X]$, $z \in \R$, defined at least for $|z|\leq 1$.  In the case of lattice-valued probability laws is often more convenient to work with the p.g.f. than the m.g.f. 

A discrete law $\mu$ belongs to the negative binomial distribution class NB$(\mathbb N)$ if its probability mass function (p.m.f.) $P_\mu$,  $P_\mu(n):=\mu(\{n\})$, and p.g.f. $\varphi_\mu$ are given respectively by
\begin{equation}
P_\mu(n)=\frac{\Gamma(\beta+n)}{\Gamma(\beta)n!}p^n (1-p)^\beta ,\quad \varphi_\mu(z)=\left(\frac{1-p}{1-p z} \right)^{\beta}, \qquad n \in \mathbb N, \, p \in (0,1), \, \beta>0.
\end{equation}
It can be shown that NB$(\mathbb N) \subset$ ID$(\R_+)$, and their associated L\'evy processes are of finite activity, enjoying various different representations, see e.g. \citet{koz+pod:09}.

Following \citet{bon:12}, Chapter 8, and \citet{ste+vh:03}, Section VI.8, 
we say that $\mu$ is a generalized negative binomial convolution, and write $\mu \in $ GNBC$(\R_+)$, if $\mu$ is obtained as a weak limit of convolutions of negative binomial laws i.e, there exists a sequence $(\mu_n)_{n \geq 0}  $, $\mu_{n}=\mu_{k_1}*\ldots *\mu_{k_n}$, $\mu_{k_i}  \in$ NB$(\mathbb N)$, such that  $\lim_{n \rightarrow \infty}\mu_n =^d\mu$. We have GNBC$(\R_+) \subset$ ID$(\mathbb R_+)$ and the p.g.f. of $\mu \in$ GNBC$(\R_+)$ writes as
\begin{equation}\label{eq:GNBC}
\varphi_\mu(z) =\exp \left( a(z-1) + \int_0^1 \log \left(\frac{1-q}{1-q z }\right) \pi(dq)\right), \quad z \in (0,1),
\end{equation}
with $a \geq 0$ and $\pi$ a positive measure supported in $[0,1]$ such that \begin{equation}\label{eq:pidef}
\int_0^{1/2} q \pi(dq)<\infty , \qquad \int_{1/2}^1 \log(1-q) \pi(dq)<\infty.
\end{equation}
It should be also noted that $\pi=h_\sharp \tau$, with $\tau$ the Thorin measure of a GGC process and $h(x)=(x+1)^{-1}$, $x \geq 0$.
An NB$(\mathbb N)$ subordinator of parameter $p$ is obtained by letting  $a=0, \pi(dq)=\beta \delta_{p}(dq), \beta>0$.

We further extend the GNBC$(\R_+)$ class to a scale-location family of lattice-valued distributions, denoted LGNBC$(\R_+)$, by considering those laws having p.g.f. of the form \begin{equation}\label{eq:LGNBC}
\varphi(z)=z^c \varphi_0(z^\alpha) , \qquad c , \alpha >0
\end{equation} with $\varphi_0$ the p.g.f. of some $\mu_0 \in$ GNBC$(\R_+)$.  Clearly LGNBC$(\R_+) \subset$ ID$(\R_+)$ and $\mu \in $ LGNBC$(\R_+)$ is characterized by the quadruplet $(c, \alpha, a, \pi)$.
An (L)GNBC$(\R_+)$ L\'evy process $Z=(Z_t)_{t \geq 0}$ is one such that $\mathcal L(Z_1) \in $ LGNBC$(\R_+)$. It is observed in \citet{bon:12}, Chapter 8, that GNBC$(\R_+)$ distributions can be characterized as a convolution of Poisson mixtures. The dynamic process counterpart, is that of an independent sum of a Poisson process $P$, whose unit time p.g.f. corresponds to the first summand in \eqref{eq:GNBC}, and the L\'evy subordination $Q_T$ of a Poisson process $Q$ to some independent Thorin subordinator $T$, with unit time  p.g.f.  given by the exponential of the integral in \eqref{eq:GNBC}. It follows that for a  LGNBC$(\R_+)$ L\'evy process $Z=(Z_t)_{t \geq 0}$ we have the representation
\begin{equation}
Z_t=^d c t + \alpha P_t +  \alpha Q_{T_t}.
\end{equation}
with $P,Q,T$ as above mutually independent. Clearly $Z$ is itself a subordinator of finite activity. When $T$ is a gamma process then $Q_T$ is simply an NB process. For details see the aforementioned Chapter 8 of  \citet{bon:12}, also providing many other noteworthy properties of GNBC distributions.

Let $Y=(X_{T_t})_{t \geq 0}$, with $X=(X_t)_{t \geq 0}$ a L\'evy process and $T=(T_t)_{t \geq 0}$ an independent LGNBC$(\R_+)$ subordinator, and consider the subordinated process $Y=(X_{T_t})_{t \geq 0}$. Letting $\mu_{X,t}=\mathcal L(X_t)$, $\mu_{Y,t}=\mathcal L(Y_t)$ and conditioning, it is plain that 
 \begin{equation}\label{eq:pfcf}
\hat \mu_{Y,1}(z)=\varphi(\hat \mu_{X,1}(z))
\end{equation}
with $\varphi$ as in \eqref{eq:LGNBC}.
In the following we argue that LGNBC$(\R_+)$ processes are the natural subordinators for $\Gamma(\R_+)$ and Bil$\Gamma(\R)$ processes, and that such subordination produces BG$_0$ Thorin processes. 

To begin with, much in the same vein of Theorem \ref{thm:ThorinSubBM1}, the Thorin measure of a subordinated gamma process to a LGNBC$(\R_+)$ process can be obtained by means of a pushforward of the Thorin measure of the subordinator.

 \begin{theoremEnd}[proof at the end,
	no link to proof]{prop}\label{prop:subgamma}
 Let $X=(X_t)_{t \geq 0}$, a $\Gamma(\alpha, \beta)$ process and $T=(T_t)_{t\geq 0}$ a {\upshape LGNBC}$(\R_+)$ subordinator independent of $X$ with triplet $(c , 1/\alpha, 0, \pi)$, and assume $\int_0^1 \pi(dq)=\vartheta >0$, $ c \, \alpha \geq \vartheta$. Then  $Y=(X_{T_t})_{t\geq 0}$ is a $ T(\R_+)$ subordinator  
 with Thorin triplet $(c',0, \tau)$ where 
\begin{align}   
c'&=\int_1^\infty \frac{\tau(dx)}{x}  \label{eq:GNBCsubLM0} \\
 \tau &= f_\sharp\pi + (c \, \alpha-\vartheta) \delta_\beta  \label{eq:GNBCsubLM}
 \end{align}
 where $f:[0,1] \rightarrow [0, \beta]$ is given by $ f(x)=\beta (1-x) $. Moreover, supp$\, \tau \subseteq [0, \beta]$ and $Y \in {\mbox {\upshape BG}}_{0}(c \, \alpha , \R_+)$.
  \end{theoremEnd}
  
\begin{proofE} In view of \eqref{eq:pfcf}, \eqref{eq:LGNBC}, \eqref{eq:GNBC},  and  using $ (\hat \mu_{X,1})^{1/\alpha}(z)=(1-iz/\beta)^{-1}=:\exp(\chi(z))$ we have
\begin{align}
\log \hat \mu_{Y,1}(z)&=   c \, \alpha \chi(z)+ \int_0^1 \log \left(\frac{1-q}{1-q (\hat \mu_{X,1}(z))^{1/\alpha} }\right) \pi(dq)  \nonumber \\ &=c \, \alpha \chi(z)+ \int_0^1 \log \left(\frac{1-q}{1-q (1- iz/\beta)^{-1} }\right) \pi(dq) \nonumber \\ &= c \, \alpha \chi(z) + \int_0^1\log \left(\frac{\beta(1-q)(1-iz/\beta)}{ \beta(1 -q)-iz  }\right) \pi(dq)  \nonumber \\ &= (c \, \alpha - \vartheta)\chi(z)+ \int_0^1\log \left(\frac{ f(q)}{ f(q)-iz  }\right)\pi(dq), \qquad z \in \R.
\end{align}
 Since $f$ continuous it holds that  supp$\,\tau \subseteq [0, 	\beta]$, so that  $\int_1^\infty x^{-1}\tau(dx) <\infty$ and   
\begin{align}\label{eq:logpitau}
\int_0^{1 \wedge \beta}|\log x|  \, \tau(dx)& \leq \int_{1-\beta^{-1}\vee 0}^1 |\log(1-q))|\pi(dq)+ \vartheta \log \beta + c \, \alpha -\vartheta \end{align} 
and the right hand is finite by assumption \eqref{eq:pidef}. This verifies \eqref{eq:thorconv0} for $\tau$ and
then \eqref{eq:GNBCsubLM0}--\eqref{eq:GNBCsubLM} follow from \eqref{eq:harexpThor}.
Finally, since $\pi$ is a finite measure and the domain of $f$ is bounded, $\tau=f_\sharp\pi$ is also finite, and  \begin{equation}
\tau_\circ(\R_+)=\tau(\R_+)=\pi(f^{-1}(\R_+))+ c \, \alpha-\vartheta=\pi([0,1])+c \, \alpha-\vartheta=c \, \alpha\end{equation} so that the last claim  follows from Definition \ref{def:GB0}.
\end{proofE}

An analogous principle holds for bilateral gamma LGNBC subordination.

\begin{theoremEnd}[proof at the end,
	no link to proof]{thm}\label{thm:subgammabil}
Let $X=(X_t)_{t \geq 0}$  be a  {\upshape Bil}$\Gamma( \alpha, \beta_+,\alpha, \beta_-)$ a bilateral gamma process and $T$ a {\upshape LGNBC}$(\R_+)$ subordinator independent of $X$ with quadruplet $(c , 1/\alpha, 0, \pi)$, satisfying $\int_0^1 \pi(dq)=\vartheta >0, c \, \alpha \geq \vartheta$. Then $Y=X_T$ is a $T(\R)$ process with Thorin triplet $(c', 0, \tau)$ given by
 \begin{align}
 c'&=\int_{\{|x|\geq 1\}}\frac{\tau(dx)}{x}  \label{eq:GNBCsubGammaLMdrift} \\
 \tau_+&= (f_\sharp\pi)^{\beta_0} + (c \, \alpha -\vartheta) \delta_{\beta_+},   \qquad  \tau_-= (f_\sharp\pi)^{-\beta_0} + (c \, \alpha-\vartheta) \delta_{\beta_-}  \label{eq:GNBCsubGammaLM}
 \end{align}
 with $\beta_0=(\beta_--\beta_+)/2$ and  $f:[0,1] \rightarrow [ |\beta_0|, ( \beta_++\beta_-)/2]$  defined as \begin{equation}
 f(x)= \frac{1}{2}\sqrt{(\beta_+ +\beta_-)^2- 4 \beta_+ \beta_- x} 
.
 \end{equation}
 Furthermore, {\upshape supp}$\,\tau  \subseteq [-\beta_-, \beta_+]$ and $Y \in \mbox{{\upshape BG}}_0(2 c \, \alpha, \R)$. 
\end{theoremEnd}

\begin{proofE}

Using the conditional argument under independence, and in view of \eqref{eq:pfcf}, \eqref{eq:GNBC} and \eqref{eq:LGNBC}  
we have, with $\exp(\chi(z)):=(\hat \mu_{X,1})^{1/\alpha}(z) = \left( (1-iz/\beta_+)(1+iz/\beta_-)\right)^{-1}$ the following equalities
\begin{align}\label{eq:pitotau}
\log \hat \mu_Y^1(z)&= c \, \alpha \chi(z) + \int_0^1 \log \left(\frac{1-q}{1-q (\hat \mu_{X,1}(z))^{1/\alpha} }\right) \pi(dq)  \nonumber \\ &=c \, \alpha \chi(z)+ \int_0^1 \log \left(\frac{1-q}{1-q (1- iz/\beta_+)^{-1}(1+ iz/\beta_-)^{-1} }\right) \pi(dq) \nonumber \\  &= c \,  \alpha \chi(z)+ \int_0^1\log \left(\frac{\beta_+ \beta_-(1-q)(1-iz/\beta_+)(1+iz/\beta_-)}{ -\beta_+ \beta_-q+(\beta_-+iz )(\beta_+-iz) }\right) \pi(dq) \nonumber \\ &= (c \, \alpha - \vartheta)\chi(z)+ \int_0^1\log \left(\frac{\beta_+ \beta_-(1-q)}{ \beta_+ \beta_-(1-q) + i z (\beta_+ -\beta_-) + z^2  }\right)\pi(dq)  \nonumber \\ 
&=(c \, \alpha - \vartheta)\chi(z) + \int_0^1\log \left(\frac{f^+(q)}{f^+(q)- iz   }\right)\pi(dq) +  \int_0^1\log \left(\frac{f^-(q)}{f^-(q) + iz   }  \right)\pi(dq)  \end{align}
with $f^{\pm}(q)= y_\pm$, and the latter are the positive solutions of \begin{equation}
\left\{ \begin{array}{cc} y_+ y_- &= (1-q)\beta_+\beta_- \\ y_+ - y_- &= \beta_+-\beta_-  
\end{array}\right.
\end{equation}
that is $ y_\pm=\mp \beta_0+\sqrt{(\beta_+ + \beta_-)^2- 4 \beta_+\beta_- q}/2.$
Now observing $s^{\beta_0}\circ f=f^+$ and  $s^{-\beta_0} \circ f=f^-$, and assuming the expression for $c'$ holds true, \eqref{eq:GNBCsubGammaLM} follows from the Thorin--Bondesson representation after operating the pushforward on the last line of \eqref{eq:pitotau}, applying \eqref{eq:tau1d}, and writing $\chi(z)=-\int_0^\infty \log(1-iz/s) \delta_{\beta_+}(ds)- \int_0^\infty \log(1 + iz/s) \delta_{\beta_-}(ds) $.
Hence, to complete the argument, observe again  $f^\pm([0,1])\subseteq[0, \beta_\pm]$,   so that {\upshape supp}$\,\tau  \subseteq [-\beta_-, \beta_+]$ by continuity of $f^\pm$.  As a consequence, $\tau(\{|x|\geq1\})<\infty$ whence $\int_{\{|x|\geq 1\}}^\infty x^{-1}\tau(dx) <\infty$ and the expression for $c'$ also follows. This partly proves \eqref{eq:thorconv0}. Moreover, for some $0 \leq C_1<1$ it holds that
\begin{align}\label{eq:ThorMeasProof}
\int_{0}^{1 \wedge \beta_+} |\log x| \, &\tau_+(dx) \leq \int_{0}^{1 \wedge \beta_+} |\log q| f^+_\sharp\pi(dq)  + c \, \alpha - \vartheta \nonumber \\ & \leq  \int_{C_1}^{1 } \left|\log\left(\beta_+-\beta_- + \sqrt{(\beta_+ +\beta_-)^2- 4 \beta_+ \beta_- q}\right)\right| \pi(dq) +	\vartheta \log 2 + c \, \alpha - \vartheta.
\end{align} When $\beta_+>\beta_-$ the integrand above is bounded in $q=C_1$ and $q=1$ and convergence follows from finiteness of $\pi$. 
However when $\beta_+<\beta_-$  it instead diverges logarithmically in the upper intergation bound. In such case we make use of the fact that expanding in Taylor  series around $ q = 1$ one obtains
\begin{equation}
 \sqrt{(\beta_+ +\beta_-)^2- 4 \beta_+ \beta_1 q}  = \beta_- -\beta_+ + (1-q)\frac{2 \beta_-\beta_+}{\beta_--\beta_+} +O((1-q)^2),
\end{equation}
hence the integrand in \eqref{eq:ThorMeasProof} is asymptotically equivalent to $-\log(1-q)-\log(2\beta_+ \beta_-/(\beta_+-\beta_-))$, which is integrable around 1 by  \eqref{eq:pidef}. The same argument applies to  $\tau_-$ with the roles of $\beta_+$ and $\beta_-$ interchanged. Lastly, in the symmetric case $\beta_+=\beta_-:=\beta$, we have $f^\pm(x)=\beta \sqrt{1-x}$ and integrability of log with respect to $f^\pm_\sharp\pi$ follows exactly as in Proposition \ref{prop:subgamma}.

 Hence we proved that $\{\tau_+, \tau_-\}$ is a Thorin family.
Finally 
\begin{align}
\tau_\circ(\R_+)&=\tau_+(\R_+)+\tau_-(\R_+)=\pi((f^+)^{-1}(\R_+))+ 
\pi((f^-)^{-1}(\R_+))\nonumber \\ &+ 2(c \, \alpha-\vartheta)=2 \pi([0,1])+ 2(c \, \alpha-\vartheta)= 2 c \, \alpha 
\end{align}
and the proof is complete.
\end{proofE}

  Theorem \ref{thm:subgammabil} can be straightforwardly extended to LGNBC$(\R_+)$ subordination of finite convolutions of bilateral gamma processes. 
  When subordinating a (bilateral) gamma processes to a LGNBC$(\R_+)$ process, an independent gamma process may arise as an additive factor even in absence of a drift in the subordinand, which is unlike the Brownian case. However, if desired, this latter process can always be removed by choosing $c \, \alpha=\vartheta$. Notice that Theorem \ref{thm:subgammabil} does not contain Proposition \ref{prop:subgamma} as a subcase; two separate statements must be made for positive and real-valued gamma subordinands.

\begin{esmp}\label{es:stochss}
\emph{Invariance of {\upshape Bil}$\Gamma(\R)$ and $\Gamma(\R_+)$ with respect to} LGNBC \emph{subordination: ``stochastic self-similarity''}. Let $X$ be a bilateral gamma process as in Theorem \ref{thm:subgammabil}. Taking as $T$ a LGNCB$(\R_+)$ subordinator of triplet $(1, 1/\alpha, 0,  \alpha \delta_p) $, 
 we obtain from Theorem \ref{thm:subgammabil} that $\tau_\pm=\alpha \delta_{f^{\pm}(p)}$, with $f^+=s^{\beta_0}\circ f$ and  $f^-=s^{-\beta_0} \circ f$. In other words, $\mathcal L(Y_1) \in  $ Bil$\Gamma(\alpha, f^+(p), \alpha, f^-(p))$. This  property of invariance of the bilateral gamma process with respect to negative binomial subordination, in the context of Proposition \ref{prop:subgamma}, was  first noticed in \citet{koz+al:06} and dubbed ``stochastic self-similarity''. Stochastic self-similarity refers to the fact that a unit scale gamma process remains gamma (but of different scale)  after subordination by some family $\{T^\lambda\}_{\lambda>1}$. Specifically,  when $X$ is a (positively-supported) $\Gamma(\alpha, \beta)$  process  and $\{T^\lambda\}_{\lambda>1}$ a family of independent LGNBC$(\R_+)$ processes with quadruplet $(1, 1/\alpha, 0,  \alpha \delta_{1-\lambda^{-1}})$,  Proposition \ref{prop:subgamma}  establishes that $Y$ is a $\Gamma( \alpha,  \lambda^{-1} \beta)$ process, confirming \citet{koz+al:06}, Proposition 4.2. The Thorin perspective reveals that stochastic self-similarity is rooted in the fact that the representing measures of both bilateral gamma and negative binomial processes are Dirac delta distributions.
\end{esmp}

\begin{esmp}
\emph{Subordination by} LGNBC\emph{s with a beta representing density}. A natural example of LGNBC$(\R_+)$ process with absolutely continuous and integrable representing measure is obtained by exploiting the beta density function. Following \citet{bon:12}, Example 8.2.2, we choose
\begin{equation}
\pi(dq)=\frac{1}{2 \pi}q^{-1/2}(\kappa-q)^{-1/2}\I_{\{q\in (0,\kappa)\}}dq, \qquad \kappa \in (0,1]
\end{equation}
which is such that $\int_0^1\pi(dq)=1/2$.
 Now set $T$ to be the LGNBC$(\R_+)$ process with quadruplet $( (2\alpha)^{-1}, 0, \alpha, \pi)$.   Assume $X$ and $X'$ to be respectively a $\Gamma(\alpha, \beta)$  subordinator and a Bil$\Gamma( \alpha, \beta_+, \alpha, \beta_-)$ process independent of $T$. Then by Proposition \ref{prop:subgamma} and Theorem \ref{thm:subgammabil}, $Y=X_T$ and $Y'=X'_T$  belong respectively to BG$_0(1/2, \R_+)$ and BG$_0(1, \R)$. Their Thorin measures $\tau$ and $\tau'$ are respectively given by
\begin{equation}
\tau(dy)=\frac{1}{2\pi}(y-\beta(1-\kappa))^{-1/2}(\beta-y)^{1/2}\I_{\{y \in ((1-\kappa)\beta, \beta)\}} dy
\end{equation}
and
\begin{align}
\tau'_\pm(dy)&=\frac{2 y \mp (\beta_+-\beta_-)  }{2\pi}\nonumber \\ &\Big((\beta_+\mp y)(\beta_-\pm y)\Big)^{-1/2}\Big(\kappa \beta_+\beta_- - (\beta_+\mp y)(\beta_-\pm y)) \Big)^{-1/2}  \I_{\{y \in ( f^{\pm}(k) , \beta_\pm)\}} dy
\end{align}
for $f^{\pm}$ as in the previous example.
Then by Proposition \ref{prop:existsmomThor}, $\mathcal L(Y_t)$ has an analytical c.f. in a neighborhood of 0 if and only if $ \kappa < 1$. Furthermore Proposition \ref{prop:existsmomThor} implies that $E[e^{p {Y_t}}]<\infty$ if and only if $p \leq (1-k)\beta$, as this latter expression coincides with the lower bound of the support of $\tau$, and in particular the critical exponential moment exists since the Thorin measures are absolutely continuous. If instead $\kappa=1$, from Proposition \ref{prop:momprop} it follows that only the moments of order $p<1/2$ are finite, and in particular $E[Y_t]=\infty$. A similar analysis reveals that the critical exponent of $Y'$ is $\min \{f^+(\kappa), f^{-}(\kappa)\}$. For $\kappa=1$ then either $f^+(\kappa)=0$ or $f^{-}(\kappa)=0$ with both holding true in the symmetric case $\beta_+=\beta_-$. In such case the Thorin density diverges as $1/y$ in 0 and no positive moment exists.

 Let us now try to explicitly compute the canonical functions of $Y$ and $Y'$. With the substitution $z= 2 y (\kappa \beta)^{-1}- 2(1-\kappa)\kappa^{-1}$  and using \citet{gra+ryz:14}, 3.364, for $Y$ we obtain 
\begin{align}\label{eq:kappabetasub}
k(x)=&\frac{1}{2\pi} \int_{(1-\kappa)\beta}^\beta \frac{e^{-xy }}{\sqrt{(y-\beta(1-\kappa) )(\beta-y)}}dy & \nonumber \\ =&\frac{e^{-\beta(1-\kappa)x}}{2 \pi}\int_{0}^2 \frac{e^{- z \beta \kappa x/2  }}{\sqrt{z  (2-z  )}}dz    = \frac{e^{-\frac{\beta}{2} (2-\kappa)x}}{2}I_0\left(-\frac{ \kappa \beta x  }{2} \right)
\end{align}
where  $I_0$ is the modified Bessel function of the first kind of order 0, which is an even function. It is derived from \citet{abr+ste:48}, 9.2.1, that $I_0(x) \sim  e^{|x|}(2 \pi x)^{-1/2}$ as $|x|\rightarrow \infty$,  whence $k(x) \sim e^{-x (1-\kappa)\beta }  x^{-1/2} (8 \pi)^{-1/2}$ as $x \rightarrow \infty$ and the former claims on the distributional moments find confirmation in the analysis of the L\'evy moments. The analogous expression to \eqref{eq:kappabetasub} for the canonical function of $Y'$ can be attained in the symmetric case $\beta_-=\beta_+=:\beta$ and for full support in $[0,1]$ of $\pi$, i.e. $\kappa=1$. We have
\begin{align}\label{eq:kappabetasubbil}
k_\pm(x)=&\frac{1}{\pi} \int_0^\beta \frac{ e^{-xy }}{\sqrt{\beta^2-y^2}}  dy  =\frac{1}{2}\left( I_0(\beta x)-L_0(\beta x) \right)
\end{align}
where $L_0$ is the modified Struve function of order 0 (\citet{gra+ryz:14}, 3.387.5, combined with the parity of $I_0$ and $L_0$). However, explicit expressions for arbitrary values of $\kappa$, $\beta_+, \beta_-$ seem out of reach, which makes a case for the study of this process class  primarily from the Thorin standpoint.
\end{esmp}

The  examples in this last section testify further the benefits of studying Thorin processes in terms of their Thorin -- as opposed to L\'evy -- triplet.

\section{Proofs}

\printProofs

\bibliographystyle{apalike}
\bibliography{Bibliography}

\begin{thebibliography}{}

\bibitem[Abramowitz and Stegun, 1948]{abr+ste:48}
Abramowitz, M. and Stegun, I.~A. (1948).
\newblock {\em Handbook of Mathematical Functions with Formulas, Graphs, and
  Mathematical Tables}, volume~55.
\newblock US Government printing office.

\bibitem[Applebaum, 2009]{app:09}
Applebaum, D. (2009).
\newblock {\em L{\'e}vy Processes and Stochastic Calculus}.
\newblock Cambridge University Press.

\bibitem[Balakrishnan, 2013]{bal:13}
Balakrishnan, N. (2013).
\newblock {\em Handbook of the Logistic Distribution}.
\newblock CRC Press.

\bibitem[Barabesi et~al., 2016]{bar+al:16a}
Barabesi, L., Cerasa, A., Perrotta, D., and Cerioli, A. (2016).
\newblock A new family of tempered distributions.
\newblock {\em Electronic Journal of Statistics}, 10:3871--3893.

\bibitem[Barndorff-Nielsen et~al., 1982]{bn+al:82}
Barndorff-Nielsen, O., Kent, J., and S{\o}rensen, M. (1982).
\newblock Normal variance-mean mixtures and $z$-distributions.
\newblock {\em International Statistical Review}, 50:145--159.

\bibitem[Barndorff-Nielsen, 1997]{bn:97}
Barndorff-Nielsen, O.~E. (1997).
\newblock Processes of normal inverse {G}aussian type.
\newblock {\em Finance and Stochastics}, 2:41--68.

\bibitem[Barndorff-Nielsen et~al., 2006]{bns+al:06}
Barndorff-Nielsen, O.~E., Maejima, M., and Sato, K. (2006).
\newblock Some classes of multivariate infinitely divisible distributions
  admitting stochastic integral representations.
\newblock {\em Bernoulli}, 12(1):1--33.

\bibitem[Barndorff-Nielsen and Shephard, 2001]{bn+she:02}
Barndorff-Nielsen, O.~E. and Shephard, N. (2001).
\newblock Non-{G}aussian {O}rnstein--{U}hlenbeck-based models and some of their
  uses in financial economics.
\newblock {\em Journal of the Royal Statistical Society: Series B},
  63:167--241.

\bibitem[Bertoin, 1996]{ber:96}
Bertoin, J. (1996).
\newblock {\em L{\'e}vy Processes}, volume 121.
\newblock Cambridge University Press.

\bibitem[Bingham et~al., 1989]{bin+al:89}
Bingham, N.~H., Goldie, C.~M., and Teugels, J.~L. (1989).
\newblock {\em Regular Variation}.
\newblock Cambridge University Press.

\bibitem[Blumenthal and Getoor, 1961]{blu+get:61}
Blumenthal, R.~M. and Getoor, R.~K. (1961).
\newblock Sample functions of stochastic processes with stationary independent
  increments.
\newblock {\em Journal of Mathematics and Mechanics}, 10:493--516.

\bibitem[Bondesson, 1992]{bon:12}
Bondesson, L. (1992).
\newblock {\em Generalized gamma convolutions and related classes of
  distributions and densities}, volume~76.
\newblock Springer Science \& Business Media.

\bibitem[Boros and Moll, 2004]{bor+mol:04}
Boros, G. and Moll, V. (2004).
\newblock {\em Irresistible Integrals: Symbolics, Analysis and Experiments in
  the Evaluation of Integrals}.
\newblock Cambridge University Press.

\bibitem[Bretagnolle, 1972]{bre:72}
Bretagnolle, J. (1972).
\newblock $p$-variation de fonctions al{\'e}atoires 2i{\`e}me partie: Processus
  {\`a} accroissements ind{\'e}pendants.
\newblock In {\em S{\'e}minaire de Probabilit{\'e}s VI Universit{\'e} de
  Strasbourg}, pages 64--71.

\bibitem[Buchmann et~al., 2017]{buc+al:17}
Buchmann, B., Kaehler, B., Maller, R., and Szimayer, A. (2017).
\newblock Multivariate subordination using generalised {G}amma convolutions
  with applications to {V}ariance {G}amma processes and option pricing.
\newblock {\em Stochastic Processes and their Applications}, 127(7):2208--2242.

\bibitem[Burridge et~al., 2014]{bur:14}
Burridge, J., Kuznetsov, A., Kwa{\'s}nicki, M., and Kyprianou, A.~E. (2014).
\newblock New families of subordinators with explicit transition probability
  semigroup.
\newblock {\em Stochastic Processes and their Applications}, 124:3480--3495.

\bibitem[Carbotti and Comi, 2020]{car+com:20}
Carbotti, A. and Comi, G.~E. (2020).
\newblock A note on {R}iemann-{L}iouville fractional {S}obolev spaces.
\newblock {\em ArXiv preprint. arXiv:2003.09515}.

\bibitem[Carr and Torricelli, 2021]{car+tor:21}
Carr, P. and Torricelli, L. (2021).
\newblock Additive logistic processes in option pricing.
\newblock {\em Finance and Stochastics}, 25:689--724.

\bibitem[Cherny and Shiryaev, 2002]{che+shi:02}
Cherny, A.~S. and Shiryaev, A.~N. (2002).
\newblock Change of time and measure for {L}\'evy processes.
\newblock Lecture Notes for the Aarhus summer school ``From L\'evy Processes to
  Semimartingales- Recent Theoretical developments and Applications to
  Finance''.

\bibitem[Cont and Tankov, 2003]{con+tan:03}
Cont, R. and Tankov, P. (2003).
\newblock {\em Financial Modelling with Jump Processes}.
\newblock Chapman and Hall/CRC Press.

\bibitem[Dickman, 1930]{dic:30}
Dickman, K. (1930).
\newblock On the frequency of numbers containing prime factors of a certain
  relative magnitude.
\newblock {\em Arkiv for matematik, astronomi och fysik}, 22(10):A--10.

\bibitem[Feller, 1971]{fel:71}
Feller, W. (1971).
\newblock {\em An Introduction to Probability Theory and its Applications}.
\newblock Wiley \& Sons.

\bibitem[Gorenflo et~al., 2020]{gor+al:14}
Gorenflo, R., Kilbas, A.~A., Mainardi, F., and Rogosin, S. (2020).
\newblock {\em Mittag-{L}effler Functions, Related Topics and Applications.
  Second Edition}.
\newblock Springer, Berlin.

\bibitem[Gradshteyn and Ryzhik, 2014]{gra+ryz:14}
Gradshteyn, I.~S. and Ryzhik, I.~M. (2014).
\newblock {\em Table of Integrals, Series, and Products}.
\newblock Academic Press.

\bibitem[Grigelionis, 2001]{gri:01}
Grigelionis, B. (2001).
\newblock Generalized $z$-distributions and related stochastic processes.
\newblock {\em Lithuanian Mathematical Journal}, 41:239--251.

\bibitem[Grigelionis, 2008]{gri:08}
Grigelionis, B. (2008).
\newblock Thorin classes of {L}{\'e}vy processes and their transforms.
\newblock {\em Lithuanian Mathematical Journal}, 48:294--315.

\bibitem[Grigelionis, 2011]{gri:11}
Grigelionis, B. (2011).
\newblock Extending the {T}horin class.
\newblock {\em Lithuanian Mathematical Journal}, 51:194--206.

\bibitem[Gupta et~al., 2024]{gup+al:24}
Gupta, N., Kumar, A., Leonenko, N., and Vaz, J. (2024).
\newblock Generalized fractional derivatives generated by {D}ickman
  subordinator and related stochastic processes.
\newblock {\em Fractional Calculus and Applied Analysis}, pages 1--37.

\bibitem[James, 2010]{jam:10}
James, L.~F. (2010).
\newblock Lamperti-type laws.
\newblock {\em The Annals of Applied Probability}, 20:1303--1340.

\bibitem[James et~al., 2008]{jam+al:08}
James, L.~F., Roynette, B., and Yor, M. (2008).
\newblock Generalized gamma convolutions, {D}irichlet means, {T}horin measures,
  with explicit examples.
\newblock {\em Probability Surveys}, 5:346--415.

\bibitem[Jurek and Vervaat, 1983]{jur+ver:83}
Jurek, Z.~J. and Vervaat, W. (1983).
\newblock An integral representation for selfdecomposable {B}anach space valued
  random variables.
\newblock {\em Zeitschrift f{\"u}r Wahrscheinlichkeitstheorie und verwandte
  Gebiete}, 62:247--262.

\bibitem[Kilbas et~al., 1993]{kil+al:93}
Kilbas, A.~A., Marichev, O., and Samko, S. (1993).
\newblock {\em Fractional Integrals and Derivatives (Theory and Applications)}.
\newblock Gordon and Breach, Switzerland.

\bibitem[Klebanov et~al., 1985]{kle+al:84}
Klebanov, L.~B., Maniya, G.~M., and Melamed, I.~A. (1985).
\newblock A problem of {Z}olotarev and analogs of infinitely divisible and
  stable distributions in a scheme for summing a random number of random
  variables.
\newblock {\em Theory of Probability and its Applications}, 29:791--794.

\bibitem[Koponen, 1995]{kop:95}
Koponen, I. (1995).
\newblock Analytic approach to the problem of convergence of truncated {L}\'evy
  flights towards the {G}aussian stochastic process.
\newblock {\em Physical Review E}, 52:1197--1199.

\bibitem[Kozubowski et~al., 2006]{koz+al:06}
Kozubowski, T.~J., Meerschaert, M.~M., and Podg\'orski, K. (2006).
\newblock Fractional {L}aplace motion.
\newblock {\em Advances in Applied Probability}, 38:451--464.

\bibitem[Kozubowski and Podg{\'o}rski, 2009]{koz+pod:09}
Kozubowski, T.~J. and Podg{\'o}rski, K. (2009).
\newblock Distributional properties of the negative binomial {L}{\'e}vy
  process.
\newblock {\em Probability and Mathematical Statistics}, 29:43--71.

\bibitem[K{\"u}chler and Tappe, 2008]{kuc+tap:08}
K{\"u}chler, U. and Tappe, S. (2008).
\newblock Bilateral gamma distributions and processes in financial mathematics.
\newblock {\em Stochastic Processes and their Applications}, 118:261--283.

\bibitem[Laverny et~al., 2021]{lev+al:21}
Laverny, O., Masiello, E., Maume-Deschamps, V., and Rulli{\`e}re, D. (2021).
\newblock Estimation of multivariate generalized gamma convolutions through
  {L}aguerre expansions.
\newblock {\em Electronic Journal of Statistics}, 15:5158--5202.

\bibitem[L{\'e}vy, 1934]{lev:34}
L{\'e}vy, P. (1934).
\newblock Sur les int{\'e}grales dont les {\'e}l{\'e}ments sont des variables
  al{\'e}atoires ind{\'e}pendantes.
\newblock {\em Annali della Scuola Normale Superiore di Pisa-Scienze Fisiche e
  Matematiche}, 3:337--366.

\bibitem[Lukacs, 1970]{luk:70}
Lukacs, E. (1970).
\newblock {\em Characteristic Functions}.
\newblock Griffin.

\bibitem[Madan et~al., 1998]{mad+al:98}
Madan, D.~B., Carr, P.~P., and Chang, E.~C. (1998).
\newblock The variance gamma process and option pricing.
\newblock {\em Review of Finance}, 2:79--105.

\bibitem[Madan and Yor, 2008]{mad+yor:08}
Madan, D.~B. and Yor, M. (2008).
\newblock Representing the {CGMY} and {M}eixner {L}{\'e}vy processes as time
  changed {B}rownian motions.
\newblock {\em Journal of Computational Finance}, 12:27.

\bibitem[Mantegna and Stanley, 1994]{man+sta:95}
Mantegna, R.~N. and Stanley, H.~E. (1994).
\newblock Stochastic process with ultraslow convergence to a {G}aussian: the
  truncated {L}\'evy flight.
\newblock {\em Physical Review Letters}, 73:2946--2949.

\bibitem[Monroe, 1972]{mon:72}
Monroe, I. (1972).
\newblock On the $\gamma$-variation of processes with stationary independent
  increments.
\newblock {\em The Annals of Mathematical Statistics}, pages 1213--1220.

\bibitem[Pakes, 2004]{pak:04}
Pakes, A.~G. (2004).
\newblock Convolution equivalence and infinite divisibility.
\newblock {\em Journal of Applied Probability}, 41(2):407--424.

\bibitem[P{\'e}rez-Abreu and Stelzer, 2012]{per+al:12}
P{\'e}rez-Abreu, V. and Stelzer, R. (2012).
\newblock A class of infinitely divisible multivariate and matrix gamma
  distributions and cone-valued generalised gamma convolutions.
\newblock {\em arXiv preprint arXiv:1201.1461}.

\bibitem[Pillai, 1990]{pil:90}
Pillai, R.~N. (1990).
\newblock On {M}ittag-{L}effler functions and related distributions.
\newblock {\em Annals of the Institute of Statistical Mathematics},
  42:157--161.

\bibitem[Pitman and Yor, 2003]{pit+yor:03}
Pitman, J. and Yor, M. (2003).
\newblock Infinitely divisible laws associated with hyperbolic functions.
\newblock {\em Canadian Journal of Mathematics}, 55(2):292--330.

\bibitem[Rosi\'nski, 2007]{ros:07}
Rosi\'nski, J. (2007).
\newblock Tempering stable processes.
\newblock {\em Stochastic Processes and their Applications}, 117:677--707.

\bibitem[Sato, 1980]{sat:80}
Sato, K. (1980).
\newblock Class {L} of multivariate distributions and its subclasses.
\newblock {\em Journal of Multivariate Analysis}, 10:207--232.

\bibitem[Sato, 1982]{sat:82}
Sato, K. (1982).
\newblock Absolute continuity of multivariate distributions of class {L}.
\newblock {\em Journal of Multivariate Analysis}, 12:89--94.

\bibitem[Sato, 1991]{sat:91}
Sato, K. (1991).
\newblock Self-similar processes with independent increments.
\newblock {\em Probability Theory and Related Fields}, 89:285--300.

\bibitem[Sato, 1999]{sat:99}
Sato, K. (1999).
\newblock {\em L\'{e}vy Processes and Infinitely Divisible Distributions}.
\newblock Cambridge University Press.

\bibitem[Sato and Yamazato, 1978]{sat+yam:78}
Sato, K. and Yamazato, M. (1978).
\newblock On distribution functions of class {L}.
\newblock {\em Zeitschrift f{\"u}r Wahrscheinlichkeitstheorie und verwandte
  Gebiete}, 43:273--308.

\bibitem[Schoutens, 2001]{sch:01}
Schoutens, W. (2001).
\newblock Meixner processes in finance.
\newblock {\em Eurandom}.

\bibitem[Sgibnev, 1990]{sgi:90}
Sgibnev, M.~S. (1990).
\newblock Asymptotics of infinitely divisible distributions on {$R$}.
\newblock {\em Siberian Mathematical Journal}, 31(1):115--119.

\bibitem[Steutel and Van~Harn, 2004]{ste+vh:03}
Steutel, F.~W. and Van~Harn, K. (2004).
\newblock {\em Infinite Divisibility of Probability Distributions on the Real
  Line}.
\newblock Dekker, New York.

\bibitem[Thorin, 1977]{tho:77}
Thorin, O. (1977).
\newblock On the infinite divisibility of the lognormal distribution.
\newblock {\em Scandinavian Actuarial Journal}, 3:121--148.

\bibitem[Thorin, 1978]{tho:78}
Thorin, O. (1978).
\newblock An extension of the notion of a generalized ${\Gamma}$-convolution.
\newblock {\em Scandinavian Actuarial Journal}, 1978:141--149.

\bibitem[Torricelli, 2024a]{tor:24}
Torricelli, L. (2024a).
\newblock On the convolution equivalence of tempered stable distributions on
  the real line.
\newblock {\em Statistics \& Probability Letters}, 207.

\bibitem[Torricelli, 2024b]{tor:24b}
Torricelli, L. (2024b).
\newblock Radially geometric stable distributions and processes.

\bibitem[Torricelli et~al., 2022]{tor+al:21}
Torricelli, L., Barabesi, L., and Cerioli, A. (2022).
\newblock Tempered positive {L}innik processes and their representations.
\newblock {\em Electronic Journal of Statistics}, 16:6313--6347.

\bibitem[Watanabe, 2008]{wat:08}
Watanabe, T. (2008).
\newblock Convolution equivalence and distributions of random sums.
\newblock {\em Probability Theory and Related Fields}, 142:367--397.

\bibitem[Wolfe, 1971]{wol:71}
Wolfe, S.~J. (1971).
\newblock On the continuity properties of {L} functions.
\newblock {\em The Annals of Mathematical Statistics}, pages 2064--2073.

\bibitem[Wolfe, 1982]{wol:82}
Wolfe, S.~J. (1982).
\newblock On a continuous analogue of the stochastic difference equation $x_n=
  \rho x_{n-1}+ b_n$.
\newblock {\em Stochastic Processes and Their Applications}, 12:301--312.

\bibitem[Yamazato, 1978]{yam:78}
Yamazato, M. (1978).
\newblock Unimodality of infinitely divisible distribution functions of class
  {L}.
\newblock {\em The Annals of Probability}, 6(4):523--531.

\bibitem[Yamazato, 1983]{yam:83}
Yamazato, M. (1983).
\newblock Absolute continuity of operator-self-decomposable distributions on
  ${R}^d$.
\newblock {\em Journal of Multivariate Analysis}, 13:550--560.

\bibitem[Zolotarev, 1963]{zol:63}
Zolotarev, V.~M. (1963).
\newblock Structure of the infinitely divisible laws of {L}-class.
\newblock {\em Lithuanian Mathematical Journal}, 3:123--140.

\end{thebibliography}

\end{document}